\documentclass[11pt]{article}

\newtheorem{theorem}{Theorem}
\newtheorem{lemma}{Lemma}
\newtheorem{proposition}{Proposition}

\newtheorem{definition}{Definition}
\newtheorem{property}{Property}

\newtheorem{conjecture}{Conjecture}

\newcommand{\lab}[1]{\label{#1}}
\newcommand{\labs}[1]{\label{#1}}
\newcommand{\labe}[1]{\label{#1}}

\newcommand{\dating}[3]{\date{{\Large {\bf #1}}\vspace{5mm}\\{\Large#3}}\label{#2}}

\usepackage{amssymb}
\def\Bbb{\mathbb}

\newcommand{\Section}[1]{\section{#1}\setcounter{equation}{0}}

\def\argmax{\mathop{\rm argmax}}

\def\Var{\mathop{\rm Var}\nolimits}

\newcommand{\cqfd}{\qquad\framebox[2.7mm]{\rule{0mm}{.7mm}}}
\newcommand{\bm}[1]{\mbox{\boldmath $#1$}}

\newcommand{\binom}[2]{\left(\begin{array}{c}#1\\#2\end{array}\right)}

\newcommand{\st}{\strut}
\newcommand{\1}{1\hskip-2.6pt{\rm l}}

\begin{document}

\title{\mbox{}\vspace{-10mm}\\{\LARGE {\bf The Brouwer Lecture 2005\vspace{3mm}\\
Statistical estimation \vspace{2mm} with model selection}}\vspace{4mm} }
\author{
{\Large Lucien Birg\'{e}\vspace{3mm}}\\{\Large Universit\'e Paris VI\vspace{1mm}}\\
{\Large Laboratoire de Probabilit\'es et Mod\`eles Al\'eatoires\vspace{1mm}}\\
{\Large U.M.R.\ C.N.R.S.\ 7599} 
}
\dating{\mbox{}}{27/04/05}{21/03/2006}      
\maketitle 
%
\begin{abstract}
The purpose of this paper is to explain the interest and importance of  (approximate) models
and model selection in Statistics. Starting from the very elementary example of histograms
we present a general notion of finite dimensional model for statistical estimation and we
explain what type of risk bounds can be expected from the use of one such model. We then give
the performance of suitable model selection procedures from a family of such models. We
illustrate our point of view by two main examples: the choice of a  partition for designing a
histogram from an $n$-sample and the problem of variable selection in the context of
Gaussian regression.
\end{abstract}

\Section{Introduction: a story of histograms\labs{H}}

\subsection{Histograms as graphical tools\labs{H1}}
Assume we are given a (large) set of real valued measurements or data $x_1,\ldots,x_n$,
corresponding to lifetimes of some human beings in a specific area, or lifetimes of some
manufactured goods, or to the annual income of families in some country, \dots . Such
measurements have a bounded range $[a,b]$ which is often known in advance (for instance
$[0,120]$ would do for lifetimes of human beings) or can be extrapolated from the data using
the extreme values. By a proper affine transformation this range can be transformed to $[0,1]$,
which we shall assume here, for the simplicity of our presentation. To represent in a
convenient, simplified, but suggestive way, this set of data, it is common to use what is
called a {\em histogram}. To design a histogram, one first chooses some finite partition
$m=\{I_0,\ldots,I_D\}$ ($D\in\Bbb{N}$) of $[0,1]$ into intervals $I_j$, generated by an
increasing sequence of endpoints $y_0=0<y_1<\ldots<y_{D+1}=1$ so that $I_j=[y_j,y_{j+1})$
for $0\le j<D$ and $I_D=[y_D,y_{D+1}]$. Then, for each $j$, one computes the number $n_j$
of observations falling in $I_j$ and one represents the data set by the piecewise constant
function $\hat{s}_{m}$ defined on $[0,1]$ by
\begin{equation}
\hat{s}_{m}(x)=\sum_{j=0}^D\frac{n_j}{n|I_j|}\1_{I_j}(x)\quad\mbox{with }
n_j=\sum_{i=1}^n\1_{I_j}(x_i)\quad\mbox{and}\quad|I_j|=y_{j+1}-y_j.
\labe{Eq-his1}
\end{equation}
Any such histogram $\hat{s}_{m}$ provides a summary of the data with three obvious
properties. It is nonnegative; its integral is equal to one ($\int_0^1\hat{s}_{m}(x)\,dx=1$) and it belongs to the $(D+1)$-dimensional linear space $V_{m}$ of
piecewise constant functions built on the partition $m$, i.e.
\begin{equation}
V_{m}=\left\{\left.t=\sum_{j=0}^Da_j\1_{I_j}\,\right|\,a_0,\ldots,a_D\in\Bbb{R}\right\}.
\labe{Eq-Vm}
\end{equation}
If the points $y_j$ are equispaced, i.e.\ all intervals $I_j$ have the same length $(D+1)^{-1}$,
the partition and the histogram are called {\em regular}. If $D\ge1$ and all intervals do not
have the same length, the partition is called {\em irregular}.

Even within this very elementary framework, some questions are in order: what
is a ``good" partition, i.e. how can one measure the quality of the representation
of the data by a histogram, and how can one choose such a good partition? One
can easily figure out that a partition with too few intervals, as compared with $n$,
will lead to an uninformative representation. Alternatively, if there are too few
data per interval the histogram may be quite erratic and meaningless. But these
are purely qualitative properties which cannot lead to a sound criterion of quality
for a partition which could be used to choose a proper one.

\subsection{Histograms as density estimators\labs{H2}}

\subsubsection{The stochastic point of view\labs{H2a}}
To go further with this analysis, we have to put the whole thing into a more mathematical
framework and a convenient one, for this type of problem, is of statistical nature. In many
situations, our data $x_i$ can be considered as successive observations of some random
phenomenon which means that $x_i=X_i(\omega)$ is the realization of a random variable $X_i$
from some probability space $(\Omega,{\cal A},\Bbb{P})$ with values in $[0,1]$ (with its Borel
$\sigma$-algebra). If we assume that the random phenomenon was stable during the
observation period and the measurements were done independently of each
other, the random variables $X_i$ can be considered as i.i.d.\ (independent and
identically distributed) with common distribution $Q$ so that
\[
\Bbb{P}[\{\omega\in\Omega\,|\,X_1(\omega)\in A_1,\ldots,X_n(\omega)\in
A_n\}]=\prod_{i=1}^nQ(A_i),
\]
for any family of Borel sets $A_1,\ldots,A_n\subset[0,1]$. Such assumptions are justified (at
least approximately) in many practical situations and $(X_1,\ldots,X_n)$ is then called an
$n$-sample from the distribution $Q$.

With this new probabilistic interpretation, $\hat{s}_{m}=\hat{s}_{m}(x,\omega)$
becomes a random function, more precisely a random element of $V_{m}$, and
(\ref{Eq-his1}) becomes
\begin{equation}
\hat{s}_{m}(x,\omega)=\sum_{j=0}^D\frac{N_j(\omega)}{n|I_j|}\1_{I_j}(x)
\quad\mbox{with }N_j(\omega)=\sum_{i=1}^n\1_{I_j}(X_i(\omega)).
\labe{Eq-his2}
\end{equation}
From now on, following the probabilistic tradition, we shall, most of the time, omit the variable
$\omega$ when dealing with random elements.

It follows from (\ref{Eq-his2}) that the random variables $N_j$ are binomial random variables
with parameters $n$ and $p_j=Q(I_j)$ and, if we assume that $Q$ has a density $s$ with 
respect to the Lebesgue measure on $[0,1]$, then $p_j=\int_{I_j}s(x)\,dx$. If $s$ also belongs to
$\Bbb{L}_2([0,1],dx)$, the piecewise constant element
$s_{m}=\sum_{j=0}^Dp_j|I_j|^{-1}\1_{I_j}$ of $\Bbb{L}_\infty([0,1],dx)$ is the orthogonal
projection of $s$ onto $V_{m}$ and
\begin{equation}
p_j=\int_{I_j}s_{m}(x)\,dx\qquad\mbox{and}\qquad
\|s-\hat{s}_m\|^2=\|s-s_{m}\|^2+\|s_{m}-\hat{s}_{m}\|^2,
\labe{Eq-his3}
\end{equation}
where $\|t\|$ denotes the $\Bbb{L}_2$-norm of $t$.

\subsubsection{Density estimators and their risk\labs{H2b}}
From a practical point of view, even if it is reasonable to assume that the variables $X_i$ are
i.i.d.\ with distribution $Q$ and density $s=dQ/dx$, this distribution is typically unknown 
and its density as well and it is often useful, in order to have an idea of the stochastic nature of
the phenomenon that produced the data, to get as much information as possible about the
unknown density $s$. For instance, comparing the shapes of lifetime densities among different
populations or their evolution with time brings much more information than merely comparing
the corresponding expected lifetimes. The very purpose of Statistics is to derive information
about the deterministic, but unknown, parameter $s$ from the stochastic, but observable, data
$X_i(\omega)$. In our problem, $\hat{s}_{m}$, which is a density, can be viewed as a random
approximation of $s$ solely based on the available information provided by the sample
$X_1,\ldots,X_n$, i.e., in statistical language, an {\em estimator} of $s$. The  distortion of the
estimated density $\hat{s}_{m}$ from the true density $s$ can be measured by the quantity
$\|s-\hat{s}_{m}\|^2$. It is clearly not the only way but this one, as seen from (\ref{Eq-his3}),
has the advantage of simplicity. Note that $\|s-\hat{s}_{m}\|^2$ is a random quantity
depending on $\omega$ as $\hat{s}_{m}$ does. In order to average out this randomness, the
statisticians often consider, as a measure of the quality of the estimator $\hat{s}_{m}$, its {\em
risk at $s$} which is the expectation of the distortion
$\|s-\hat{s}_{m}\|^2$ given by
\[
R(\hat{s}_{m},s)=\Bbb{E}_s\left[\left\|s-\hat{s}_{m}\right\|^2\right]=
\int\left\|s-\hat{s}_{m}(\omega)\right\|^2d\Bbb{P}_s(\omega).
\]
Here $\Bbb{P}_s$ and $\Bbb{E}_s$ respectively denote the probability and the expectation of
functions of $X_1,\ldots,X_n$ when these variables are i.i.d.\ with density $s$. Of course, due
to randomness, $R(\hat{s}_{m},s)$ does not provide any information on the actual distortion
$\left\|s-\hat{s}_{m}(\omega)\right\|^2$ in our experiment. But, by the law of large
numbers, it provides a good approximation of the average distorsion one would get if one
iterated many times the procedure of drawing a sample $X_1,\ldots,X_n$ and building the
corresponding histogram. The importance of the risk, as a measure of the quality of the
estimator $\hat{s}_{m}$ also derives from Markov Inequality which implies that, for any $z>0$,
\begin{equation}
\Bbb{P}_s\left[\left\|s-\hat{s}_{m}\right\|\ge \sqrt{zR(\hat{s}_{m},s)}\right]\le
z^{-1}.
\labe{Eq-his4}
\end{equation}
Hence, with a guaranteed probability $1-z^{-1}$, the distance between $s$ and its estimator is
bounded by $\sqrt{zR(\hat{s}_{m},s)}$. When $z$ is large, there are only two cases: either
we were very unlucky and an event of probability not larger than $z^{-1}$ occurred, or we were
not and $\left\|s-\hat{s}_{m}\right\|\le \sqrt{zR(\hat{s}_{m},s)}$. Of course, there
is no way to know which of the two cases occured, but this is the rule in Statistics: there is
always some uncertainty in our conclusions.

\subsubsection{Risk bounds for histograms\labs{H2c}}
In any case, (\ref{Eq-his4}) shows that the risk can be viewed as a good indicator of the
performance of an estimator. Moreover, it follows from (\ref{Eq-his3}) that it can be written as
\begin{equation}
R(\hat{s}_{m},s)=\|s-s_{m}\|^2+
\Bbb{E}_s\left[\left\|s_{m}-\hat{s}_{m}\right\|^2\right].
\labe{Eq-his8}
\end{equation}
With this special choice of distortion, the risk can be decomposed into the sum of two terms. The
first one has nothing to do with the stochastic nature of the observations but simply measures
the quality of approximation of $s$ by the linear space $V_{m}$ since it is the square of the
distance from $s$ to $V_{m}$. It only depends on the partition and the true unknown density
$s$, not on the observations. 

The second term in the risk, which is due to the stochastic nature of the observations, hence of
$\hat{s}_m$, can be bounded in the following way, since $N_j$ is a binomial random variable
with parameters $n$ and $p_j$ and both $s_{m}$ and $\hat{s}_{m}$ are constant on each interval
$I_j$:
\begin{eqnarray}
\Bbb{E}_s\left[\left\|s_{m}-\hat{s}_{m}\right\|^2\right]&=&\sum_{j=0}^{D}
\Bbb{E}_s\left[\int_{I_j}\left(s_{m}(x)-\hat{s}_{m}(x)\right)^2\,dx\right]
\nonumber\\&=&\sum_{j=0}^{D}\Bbb{E}_s\left[|I_j|
\left(\frac{p_j}{|I_j|}-\frac{N_j}{n|I_j|}\right)^2\right]\nonumber\\&=&
\sum_{j=0}^{D}\frac{1}{n^2|I_j|}\Var(N_j)\;\;=\;\;\
\frac{1}{n}\sum_{j=0}^{D}\frac{p_j(1-p_j)}{|I_j|}.
\labe{Eq-his0}
\end{eqnarray}
This quantity is easy to bound in the special case of a regular partition since then
$|I_j|=(D+1)^{-1}$ and we get, using the concavity of the function $x\mapsto x(1-x)$,
\begin{eqnarray}
\Bbb{E}_s\left[\left\|s_{m}-\hat{s}_{m}\right\|^2\right]&=&\frac{(D+1)^2}{n}
\sum_{j=0}^{D}\frac{p_j(1-p_j)}{D+1}\nonumber\\&\le&\frac{(D+1)^2}{n}
\frac{\sum_{j=0}^Dp_j}{D+1}\left(1-\frac{\sum_{j=0}^Dp_j}{D+1}\right)\;\;=\;\;\frac{D}{n}.
\labe{Eq-his9}
\end{eqnarray}
Note that $D=0$ corresponds to the degenerate partition $m_0=\{[0,1]\}$ for which
$s_{m_0}=\1_{[0,1]}$ which is the density of the uniform distribution on $[0,1]$,
independently of $s$. Then $\hat{s}_{m_0}=s_{m_0}=\1_{[0,1]}$ and
$R(\hat{s}_{m_0},s)=\|s-\1_{[0,1]}\|^2$.

For general irregular partitions we derive from (\ref{Eq-his3}) that 
$p_j\le|I_j|\|s_{m}\|_\infty$, hence, by (\ref{Eq-his0}),
\begin{equation}
\Bbb{E}_s\left[\left\|s_{m}-\hat{s}_{m}\right\|^2\right]\le\frac{\|s_{m}\|_\infty}{n}\sum_{j=0}^{D}(1-p_j)=\frac{D\|s_{m}\|_\infty}{n}.
\labe{Eq-his5}
\end{equation}
There is actually little space for improvement in (\ref{Eq-his5}) as shown by the following
example. Define the partition $m$ by $I_j=[\alpha j,\alpha(j+1))$ for $0\le j<D$ and
$I_D=[\alpha D,1]$ with $0<\alpha<D^{-1}$. Set $s=s_{m}=(\alpha D)^{-1}
\left(1-\1_{I_D}\right)$. Then $p_j=D^{-1}$ for  $0\le j<D$ and, by (\ref{Eq-his0}),
\[
\Bbb{E}_s\left[\left\|s_{m}-\hat{s}_{m}\right\|^2\right]=\frac{D-1}{\alpha Dn} =\frac{(D-1)\|s_{m}\|_\infty}{n}.
\]
If we make the extra assumption that $s$ belongs to $\Bbb{L}_\infty([0,1],dx)$, then
$\|s_{m}\|_\infty\le\|s\|_\infty$ and (\ref{Eq-his5}) becomes $\Bbb{E}_s
\left[\left\|s_{m}-\hat{s}_{m}\right\|^2\right]\le \|s\|_\infty n^{-1}D$. This bound is  also
valid for regular partitions but always worse than (\ref{Eq-his9}) since $\|s\|_\infty\ge1$ for
all densities with respect to Lebesgue measure on $[0,1]$ and strictly worse if $s$ is not the
uniform density.  Finally, by (\ref{Eq-his8}),
\begin{equation}
R(\hat{s}_{m},s)\le\|s-s_{m}\|^2+\|s\|_\infty n^{-1}D.
\labe{Eq-his7}
\end{equation}
As we shall see later the rather unpleasant presence of
the unknown and possibly unbounded $\|s\|_\infty$ factor in the second term is due to the
way we measure the distance between densities, i.e.\ through the $\Bbb{L}_2$-norm. 

\subsection{A first approach to model selection\labs{H3}}

\subsubsection{An alternative interpretation of histograms\labs{H3a}}
The decomposition (\ref{Eq-his3}) suggests another interpretation for the construction of
$\hat{s}_{m}$. What do we do here? Since $s$ is possibly a complicated object, we replace
it by a much simpler one $s_{m}$ and estimate it by $\hat{s}_{m}$. Note that $s_{m}$ is
unknown,  as $s$ is, and what is available  to the statistician is the partition $m$, the
corresponding linear space $V_{m}$ and, consequently, the set $S_{m}$ of all densities
belonging to $V_{m}$, i.e.
\begin{equation}
S_{m}=\left\{\left.t=\sum_{j=0}^Da_j\1_{I_j}(x)\,\right|\,a_0,\ldots,a_D\in
\Bbb{R}_+\quad\mbox{and}\quad\sum_{j=0}^Da_j|I_j|=1\right\}.
\labe{Eq-Sm}
\end{equation}
It is a convex subset of some $D$-dimensional linear space and $s_{m}$ is given by
$\|s-s_{m}\|=\inf_{t\in S_{m}}\|s-t\|$. It is the best approximation of $s$ in
$S_{m}$. As to $\hat{s}_{m}$ it only depends on the set $S_{m}$ and the
observations in the following way, as can easily be checked:
\[
\hat{s}_{m}=\argmax_{t\in S_{m}}\sum_{i=1}^n\log(t(X_i)),
\]
which means that it maximizes the so-called {\em likelihood function}
$t\mapsto\prod_{i=1}^nt(X_i)$ for $t\in S_{m}$, the likelihood at $t$ being the joint density of
the sample computed at the observations. The estimator $\hat{s}_{m}$ is called the {\em
maximum likelihood estimator} (m.l.e.\ for short) with respect to $S_{m}$. Note that, if
$s=s_{m}$ actually belongs to $S_{m}$, the m.l.e.\ converges in probability to $s$ at rate at least
as fast as $n^{-1/2}$ when $n$ goes to infinity since then, by  (\ref{Eq-his4}), (\ref{Eq-his8}) and
(\ref{Eq-his5}),
\[
\Bbb{P}_s\left[\left\|s-\hat{s}_{m}\right\|\ge 
n^{-1/2}\sqrt{z\|s_{m}\|_\infty D}\right]\le z^{-1}.
\]
The m.l.e.\ therefore appears to be a suitable estimator to use if the model $S_{m}$ is
correct, i.e.\ if $s\in S_{m}$. When we use the histogram estimator $\hat{s}_{m}$,
we just do as if $s$ did belong to $S_{m}$, using $S_{m}$ as an {\em approximate
model} for $s$. The resulting risk is then the sum of two terms, an approximation error equal
to the square of the distance from $s$ to $S_{m}$ and due to the fact that $s$ does not in
general belong to the model $S_{m}$, and an estimation term
$\Bbb{E}_s\left[\left\|s_{m}-\hat{s}_{m}\right\|^2\right]$ which is the risk
corresponding to the estimation within  the model when $s=s_{m}$ since
$\left\|s_{m}-\hat{s}_{m}\right\|^2$ has the same expectation when the
observations are i.i.d.\ with density $s$ or $s_{m}$. 

\subsubsection{Model selection and oracles\labs{H3b}}
Let us denote by $m_D$ the regular partition with $D+1$ pieces and set
$S_D=S_{m_D}$, $\hat{s}_D=\hat{s}_{m_D}$ and $s_D=s_{m_D}$, for
simplicity. It follows from (\ref{Eq-his8}) and (\ref{Eq-his9}) that
\begin{equation}
R(\hat{s}_D,s)\le\|s-s_D\|^2+n^{-1}D.
\labe{Eq-his6}
\end{equation}
From the approximation point of view, a good partition should lead to a small value of
$\|s-s_D\|$ which typically requires a partition into many intervals, hence a large
value of $D$, while the estimation point of view requires a model $S_D$ defined by few
parameters, hence a small value of $D$. Obviously, these requirements are contradictory and
one should look for a compromise between them in order to minimize the right-hand side of
(\ref{Eq-his6}). Unfortunately, the value $D_{opt}$ which satisfies
\[
\|s-s_{D_{opt}}\|^2+n^{-1}D_{opt}=\inf_{D\in\Bbb{N}}\left\{\|s-s_D\|^2+n^{-1}D\right\}
\]
cannot be computed since it depends on the unknown density $s$ via the approximation term
$\|s-s_D\|$ and is not accessible to the statistician. This is why the random variable
$\hat{s}_{D_{opt}}$ based on the partition $m_{D_{opt}}$ is called an ``oracle". It is not an
estimator because it makes use of the number $D_{opt}$ which is unknown to the statistician.
The problem of model selection is to find a genuine estimator, solely based on the data, that
mimics an oracle, i.e.\ to use the data $X_1,\ldots,X_n$ to select a number
$\widehat{D}(X_1,\ldots,X_n)$ such that the resulting histogram $\tilde{s}=
\hat{s}_{\widehat{D}}$ has a performance which is comparable to that of the oracle:
\[
R(\tilde{s},s)\le C\left[\|s-s_{D_{opt}}\|^2+n^{-1}D_{opt}\right],
\]
where $C$ is a constant that neither depends on the unknown density $s$ nor on $n$.

\subsubsection{An illustrative example\labs{H3c}}
Still working with the regular partitions $m_D$, let us now assume that the unknown density
$s$ satisfies some H\"olderian continuity condition, 
\begin{equation}
|s(x)-s(y)|\le L|x-y|^\beta,\quad L>0,\quad 0<\beta\le1\quad\mbox{for all }x,y\in[0,1].
\labe{Eq-Hold}
\end{equation}
If $0\le j\le D$ and $x\in I_j$, then $s_D(x)=s(y)$ for some $y\in I_j$, hence $|s(x)-s_D(x)|\le
L(D+1)^{-\beta}$, from which we derive that $\|s-s_D\|^2\le\|s-s_D\|_\infty^2\le
L^2(D+1)^{-2\beta}$. Therefore (\ref{Eq-his6}) implies that $R(\hat{s}_D,s)\le n^{-1}D+
L^2(D+1)^{-2\beta}$. Since the minimum of the function $x\mapsto n^{-1}x+L^2x^{-2\beta}$
is obtained for $x=\left(2\beta nL^2\right)^{1/(2\beta+1)}$, we choose $D$ so that $D+1$ is
the smallest  integer $\ge\left(nL^2\right) ^{1/(2\beta+1)}$. If $nL^2\le1$, this leads to
$D=0$ and $R(\hat{s}_0,s)\le L^2\le n^{-1}$. Otherwise, and this necessarily happens
for large enough $n$, $1\le D<\left(nL^2\right)^{1/(2\beta+1)}$, hence
$R(\hat{s}_D,s)\le2\left(Ln^{-\beta}\right) ^{2/(2\beta+1)}$. Finally, in any case,
\[
R(\hat{s}_D,s)\le\max\left\{2\left(Ln^{-\beta}\right)^{2/(2\beta+1)};n^{-1}\right\}.
\]
Unfortunately, we can only get a risk bound of this form if we fix $D$ as a function of $L$ and
$\beta$, as indicated above.  Typically, $L$ and $\beta$ are also unknown so that we do not
know how to choose $D$ and cannot get the right risk bound. The situation is even more
complicated since, for a given $s$, there are many different pairs $L,\beta$ that satisfy
(\ref{Eq-Hold}), leading to different values of $D$ and risk bounds. Of course, one would like
to choose the optimal one which means choosing the value of $D$ that minimizes the
right-hand side of (\ref{Eq-his6}).

\subsection{A brief summary of this paper\labs{H4}}
The study of histograms as density estimators shows us that a convenient method to
estimate a complicated object as a density $s$ on $[0,1]$ works as follows: choose an
approximate model $S_{m}$ for $s$ involving only a limited number of unknown parameters
and then do as if the model were correct, i.e.\ if  $s\in S_{m}$, using an estimator $\hat{s}_{m}$
which is a good estimator when the model is actually correct. The resulting risk is the sum of an
approximation term which measures the quality of approximation of $s$ by the model and an
estimation term which is roughly proportional to the number of parameters needed to describe
an element of the model, reflecting its complexity. As a consequence, a good model should be
simple (described by few parameters) and accurate (close to the true density $s$). Unfortunately,
because of the second requirement, a theoretical choice of a good model should be based on the
knowledge of $s$. Given a family of possible models, a major problem is therefore to understand
to what extent one can guess from the data which model in the family is appropriate.

The remainder of this paper is devoted to giving some hints to justify and understand the
various steps needed to formally develop the previous arguments. The next section will present
the classical parametric theory of estimation which assumes that one works with the correct
model and that this model satisfies some specific regularity conditions. Under such conditions
the m.l.e.\ enjoys some good asymptotic properties that we shall recall, but this classical theory
does not handle the case of approximate models or infinite dimensional parameters. It has
therefore been extended in the recent years in many directions to (partly) cover such situations.
We shall present here one such generalization that attempts to solve (at least theoretically) most
of the difficulties connected with the classical theory. In Section~\ref{D}, we shall depart from
the classical theory, assuming only an approximate model and checking on some examples that
the results we got for histograms essentially extend to these cases with a risk  bounded by an
approximation term plus an estimation term which again leads to the problem of selecting a
good model. Section~\ref{A} is devoted to a more general approach to estimation based on an
approximate model with {\em finite dimension} for a suitably defined and purely metric
notion of dimension. We show here that some specific estimators (sometimes discretized
versions of the m.l.e., sometimes more complicated ones) do lead to risk bounds of the
required form: an approximation term plus an estimation term which is proportional to the
{\em dimension} (when suitably defined) of the model. In the last section, we explain how to
handle many such approximate models with finite dimensions simultaneously. Ideally, we
would like to choose, using only the data, the best model in the family, i.e.\ the one with the
smallest risk. This is unfortunately not possible, but we shall explain to what extent one can
approximate this ideal risk.

\Section{Some historical considerations\labs{C}}

\subsection{The classical parametric point of view\labs{C1}}
To be specific, let us assume again that our observations $X_1,\ldots,X_n$ are i.i.d.\ random
variables with an unknown density $s$ with respect to some reference measure $\nu$
defined on the underlying measurable set $(E,{\cal E})$ (not necessarily the Lebesgue
measure on $[0,1]$) so that the joint distribution $P_s$ of the observations on $E^n$ is
given by 
\[
\frac{dP_s}{d\nu^{\otimes n}}(x_1,\ldots,x_n)=\prod_{i=1}^ns(x_i).
\]
In the sequel, we shall call the problem of estimating the unknown density $s$ from the i.i.d.\
sample $X_1,\ldots,X_n$ the {\em density estimation problem} or the {\em i.i.d.\
framework}.

The classical parametric approach to density estimation that developed after milestone papers
by Fisher (1921 and 1925) up to the sixties and is still quite popular nowadays is somewhat
different from what we described before. It typically  assumes a {\em parametric model}
$\overline{S}$ for $s$, which means that the true unknown density $s$ of our observations
belongs to some particular set $\overline{S}=\{t_\theta\,|\,\theta\in\Theta\}$ of densities
parametrized by some subset $\Theta$ of a Euclidean space $\Bbb{R}^k$. Then
$s=t_{\theta_0}$ for some particular $\theta_0\in\Theta$ which is called the {\em true
parameter value}. One assumes moreover that the mapping $\theta\mapsto t_\theta$ from
$\Theta$ to $\overline{S}$ is smooth (in a suitable sense) and one-to-one, so that estimating $s$
is equivalent to estimating the parameter $\theta_0$. An estimator
$\hat{\theta}_n(X_1,\ldots,X_n)$ of $\theta_0$ is then defined via a measurable mapping
$\hat{\theta}_n$ from $E^n$ to $\Theta$ (with its Borel $\sigma$-algebra) and its quadratic
risk is given by
\[
R(\hat{\theta}_n,\theta_0)=\Bbb{E}_s\left[\|\hat{\theta}_n-\theta_0\|^2\right],
\]
where $\|\cdot\|$ now denotes the Euclidian norm in $\Bbb{R}^k$. Typical examples of
parametric models for densities on the real line are given by 

i) the Gaussian densities ${\cal N}\left(\mu,\sigma^2\right)$ with
$\theta=\left(\mu,\sigma^2\right)$ and $\Theta=\Bbb{R}\times(0,+\infty)$ given by
\[ 
t_\theta(x)=\frac{1}{\sqrt{2\pi\sigma^2}}\exp\left[-\frac{1}{2\sigma^2}(x-\mu)^2\right];
\] 

ii) the gamma densities $\Gamma(v,\lambda)$ with $\theta=(v,\lambda)$ and
$\Theta=(0,+\infty)^2$ given by
\[ 
t_\theta(x)=[\Gamma(v)]^{-1}\lambda^vx^{v-1}\exp[-\lambda x];
\] 

iii) the uniform density on the interval $[\theta,\theta+1]$ given by
$t(\theta)=\1_{[\theta,\theta+1]}$ with $\theta\in\Bbb{R}$.

\subsection{The maximum likelihood method\labs{C2}}

\subsubsection{Consistency and asymptotic normality of the parametric m.l.e.\labs{C2a}}
Fisher's approach to parametric estimation is mainly connected with the method of maximum
likelihood. We recall from Section~\ref{H3a} that the likelihood function on $\Theta$ is
given by $\theta\mapsto\prod_{i=1}^nt_\theta(X_i)$ and a maximum likelihood estimator
$\hat{\theta}_n$ is any maximizer of this function or equivalently of the
{\em log-likelihood function}
\[
L(\theta)=\sum_{i=1}^n\log\left(\st t_\theta(X_i)\right).
\]
For Gaussian densities, the maximum likelihood estimator $\hat{\theta}_n=
\left(\hat{\mu}_n,\hat{\sigma}_n^2\right)$ is unique and given by $\hat{\mu}_n=
n^{-1}\sum_{i=1}^nX_i$ and $\hat{\sigma}_n^2=n^{-1}\sum_{i=1}^n
\left(X_i-\hat{\mu}_n\right)^2$. Moreover $\hat{\theta}_n$ converges in probability to the
true parameter $\theta_0$ when $n$ goes to infinity. We say that $\hat{\theta}_n$ is {\em
consistent}. Unfortunately, this situation is not general. The study of our second and third
examples show that explicit computation of the m.l.e.\ is not always possible (gamma
densities) or the m.l.e.\ may not be unique (uniform densities). One can also find examples of
{\em inconsistency} of the m.l.e., but, as shown by Wald (1949), it can be proved that, under
suitably strong assumptions, any sequence of maximum likelihood estimators is consistent.

If the mapping $\theta\mapsto l_\theta(x)=\log\left(\st t_\theta(x)\right)$ satisfies
suitable differentiability assumptions, the parametric model is called {\em regular}. This is
the case for the Gaussian and gamma densities, not for the uniform. If the model is regular
and the m.l.e.\ is consistent we can expand the derivative of the function $L$ in a vicinity of
$\theta_0$ when it is an inner point of $\Theta$. Restricting ourselves, for simplicity, to the
case $\Theta\subset\Bbb{R}$, we get
\[
L'(\theta)=L'(\theta_0)+(\theta-\theta_0)L''(\theta_0)+(1/2)(\theta-\theta_0)^2
L'''(\theta')
\]
and since $\hat{\theta}_n$ is a maximizer for $L$,
\[
L'(\hat{\theta}_n)=0=L'(\theta_0)+(\hat{\theta}_n-\theta_0)
L''(\theta_0)+(1/2)(\hat{\theta}_n-\theta_0)^2L'''(\theta'_n),
\]
for some sequence $(\theta'_n)$ converging to $\theta_0$ in probability as $\hat{\theta}_n$
does. Equivalently, setting $\delta_n=\sqrt{n}\left(\hat{\theta}_n-\theta_0\right)$,
\begin{equation}
0=\frac{1}{\sqrt{n}}\sum_{i=1}^nl'_{\theta_0}(X_i)+\left[\frac{1}{n}\sum_{i=1}^n
l''_{\theta_0}(X_i)\right]\delta_n+\frac{\hat{\theta}_n-\theta_0}{2}\left[\frac{1}{n}
\sum_{i=1}^nl'''_{\theta'_n}(X_i)\right]\delta_n.
\labe{Eq-like1}
\end{equation}
Since $\int t_\theta(x)\,d\nu(x)=1$ for all $\theta$, it follows from the regularity
assumptions that
\[
\Bbb{E}_s\left[l'_{\theta_0}(X_i)\right]=\int[t'_{\theta_0}(x)/t_{\theta_0}(x)]t_{\theta_0}(x)
\,d\nu(x)=\int t'_{\theta_0}(x)\,d\nu(x)=0
\]
and
\begin{eqnarray*}
\Bbb{E}_s\left[l''_{\theta_0}(X_i)\right]&=&\int[t''_{\theta_0}(x)/t_{\theta_0}(x)]t_{\theta_0}(x)
\,d\nu(x)-\int[t'_{\theta_0}(x)/t_{\theta_0}(x)]^2t_{\theta_0}(x)\,d\nu(x)\\&=&0-
\int\left([t'_{\theta_0}(x)]^2/t_{\theta_0}(x)\right)\,d\nu(x)\;\;=\;\;-I(\theta_0),
\end{eqnarray*}
where the last equality defines the {\em Fisher Information} $I(\theta_0)$. Moreover
\[
\Var\left[l'_{\theta_0}(X_i)\right]=\Bbb{E}_s\left[\left(l'_{\theta_0}(X_i)\right)^2\right]
=\int[t'_{\theta_0}(x)/t_{\theta_0}(x)]^2t_{\theta_0}(x)\,d\nu(x)=I(\theta_0).
\]
It then follows from the law of large numbers that
\[
\frac{1}{n}\sum_{i=1}^nl''_{\theta_0}(X_i)\stackrel{P}\rightarrow
\Bbb{E}_s\left[l''_{\theta_0}(X_i)\right]=-I(\theta_0)
\]
and from the central limit theorem that
\[
\frac{1}{\sqrt{n}}\sum_{i=1}^nl'_{\theta_0}(X_i)\leadsto{\cal N}(0,I(\theta_0)),
\]
where $\stackrel{P}\rightarrow$ and $\leadsto$ denote respectively the convergences in
probability and in distribution. The regularity assumptions also ensure that
$n^{-1}\sum_{i=1}^nl'''_{\theta'_n}(X_i)$ is asymptotically bounded so that the third term in
(\ref{Eq-like1}) is asymptotically negligible as compared to the other two. We finally deduce
from (\ref{Eq-like1}) that 
\begin{equation}
\delta_n=\sqrt{n}(\hat{\theta}_n-\theta_0)\leadsto{\cal N}\left(0,[I(\theta_0)]^{-1}\right).
\labe{Eq-like2}
\end{equation}
This is the so-called asymptotic normality and efficiency of the maximum likelihood
estimator and a formal proof of this result can be found in Cram\'er (1946, Section~33.3). It
can also be proved that the asymptotic variance $[I(\theta_0)]^{-1}$ of $\delta_n$ is, in various
senses, optimal, as shown by Le Cam (1953) and Hajek (1970 and 1972). Much less restrictive
conditions of regularity which still imply the asymptotic normality and efficiency of the
m.l.e.\ have been given by Le Cam (1970) --- see also Theorem~12.3 in van der Vaart (2002)
---. A good account of the theory can be found in Ibragimov and Has'minskii (1981). A more
recent point of view on the theory of regularity and the m.l.e., based on empirical process
theory, is to be found in van der Vaart (1998).

\subsubsection{A more general point of view on the maximum likelihood method
\labs{C2b}} 
The limitations of the classical parametric theory of maximum likelihood have been
recognized for a long time. We already mentioned problems of inconsistency. Examples and
further references can be found in Le Cam (1990). Moreover, although it is widely believed
among non-specialists that (\ref{Eq-like2}) typically holds, this is definitely not true, even
under consistency. For instance, if $t_\theta=\theta^{-1}\1_{[0,\theta]}$ is the uniform
density on $[0,\theta]$ and $\Theta=(0,+\infty)$, the m.l.e.\ satisfies
$n(\theta_0-\hat{\theta}_n)\leadsto\Gamma(1,\theta_0)$. Additional examples
can be found in Ibragimov and Has'minskii (1981, Chapters~5 and 6) showing that neither the
rate $\sqrt{n}$ nor the limiting normal distribution are general.

Another drawback of the classical point of view on maximum likelihood estimation is its
purely asymptotic nature. Not only does it require specific assumptions and can fail under
small departures from these assumptions but it tells us nothing about the real performances
of the  m.l.e.\ for a given (even large) number $n_0$ of observations, just as the central limit
theorem does. Suppose that our observations $X_1,\ldots,X_n$ are i.i.d.\ Bernoulli variables
taking only the values 0 and 1 with respective probabilities $1-\theta_0$ and $\theta_0$ and
$\Theta=[0,1]$. Then $\hat{\theta}_n=n^{-1}\sum_{i=1}^nX_i$ and, if $0<\theta_0<1$,
$\delta_n=\sqrt{n}(\hat{\theta}_n-\theta_0)\leadsto{\cal N}(0,\theta_0[1-\theta_0])$ as
expected. But it is well-known that if $n=1000$ and $0<\theta_0\le0.002$, the distribution of
$\hat{\theta}_n$ looks rather like a Poisson distribution with parameter $1000\theta_0$
than like a normal ${\cal N} (\theta_0,n^{-1}\theta_0[1-\theta_0])$ as predicted by the
asymptotic theory. A discussion about the relevance of the asymptotic point of view for
practical purposes can be found in Le Cam and Yang (2000, Section~7.1).

A further limitation of the classical m.l.e.\ theory is the fact that the assumed parametric
model is true, i.e.\ the unknown distribution of the observations has a density $s$ with respect
to $\nu$ which is of the form $t_{\theta_0}$ for some $\theta_0\in\Theta$. If this
assumption is violated, even slightly, the whole theory fails as can be seen from the following
example. We assume a Gaussian distribution $P_\theta={\cal N}(\theta,1)$ with density
$t_\theta$ with respect to the Lebesgue measure and $\Theta=\Bbb{R}$ but the observations
actually follow the distribution $Q=(99P_0+P_{300})/100$. It is actually rather close to the
$P_0$ distribution, which belongs to the model, in the sense that, for any measurable set $A$,
$|Q(A)-P_0(A)|\le1/100$. Nevertheless, the m.l.e.\ $n^{-1}\sum_{i=1}^nX_i$ converges to
$3$ so that the estimated distribution based on the wrong model will be close to $P_3$,
hence quite different from the true distribution which is close to $P_0$. 

For all these reasons, the classical approach to maximum likelihood estimation has been
substantially generalized in the recent years. Nonparametric and semiparametric maximum
likelihood allows to deal with families of distributions $P_s$ where $s$ belongs to some
infinite-dimensional set, while sieved m.l.e.\ involves situations where the true distribution
does not belong to the model. Both extensions lead to truely nonasymptotic results. Among the
many papers dealing with such extensions, let us mention here Grenander (1981), Silverman
(1982), Wahba (1990), Groeneboom and Wellner (1992), van de Geer (1993, 1995 and 2000),
Birg\'e and  Massart (1993 and 1998), Shen and Wong (1994), Wong and Shen (1995), van der
Vaart and Wellner (1996), Barron, Birg\'e and Massart (1999) and Massart (2006).  Let us now 
explain what are the novelties brought by some of these extentions.

\Section{An alternative point of view\labs{D}}

\subsection{Nonparametric density estimation\labs{D1}}
The assumption that the unknown density $s$ of the observations belongs to a parametric
model, i.e.\ a smooth image of some subset of a Euclidean space, appears to be definitely too
strong and unsatisfactory in many situations. Let us give here two illustrations. If we assume
that $s$ belongs to the set ${\cal S}_1$ of Lipschitz densities on $[0,1]$ (i.e.\ $s$ satisfies
$|s(x)-s(y)|\le|x-y|$), one cannot represent ${\cal S}_1$ in a smooth way by a finite number
of real parameters. The same holds if we simply assume that $s\in{\cal S}_2$, the set of all
densities in $\Bbb{L}_2([0,1],dx)$. In this case, given some orthonormal basis
$(\varphi_j)_{j\ge1}$ of $\Bbb{L}_2([0,1],dx)$, there exists a natural parametrization of
${\cal S}_2$ by $\bm{l}_2(\Bbb{N}^\star)$ ($\Bbb{N}^\star=\Bbb{N}\setminus\{0\}$) via
the coordinates, but it is definitely not finite-dimensional. These two problems are examples of
{\em nonparametric density estimation problems}.

\subsubsection{Projection estimators\labs{D1a}}
In order to solve the second estimation problem, Cencov (1962) proposed a general class of
estimators called {\em projection estimators}. The idea is to estimate the coefficients $s_j$ of
$s$ in the orthonormal expansion $s=\sum_{j=1}^{+\infty}s_j\varphi_j$ using estimators
$\hat{s}_j$ chosen in such a way that $\sum_{j=1}^{+\infty}\hat{s}_j^2<+\infty$ a.s.\ so
that $\hat{s}=\sum_{j=1}^{+\infty}\hat{s}_j\varphi_j$ belongs to $\Bbb{L}_2([0,1],dx)$
a.s.. Since $s_j=\int_0^1s(x)\varphi_j(x)\,dx=\Bbb{E}_s[\varphi_j(X_i)]$, a natural estimator
for $s_j$ is $\overline{\varphi}_j=n^{-1}\sum_{i=1}^n\varphi_j(X_i)$. Indeed 
\begin{equation}
\Bbb{E}_s\left[\overline{\varphi}_j\right]=s_j\quad\mbox{and}\quad
\Var\left(\overline{\varphi}_j\right)=n^{-1}\Var(\varphi_j(X_1))\le
n^{-1}\int_0^1\varphi_j^2(x)s(x)\,dx.
\labe{Eq-proj0}
\end{equation}
Assuming, for simplicity, that we take for $(\varphi_j)_{j\ge0}$ the trigonometric basis which
is bounded by $\sqrt{2}$, we derive that $\Var\left(\overline{\varphi}_j\right)\le2/n$. 
We cannot use $\sum_{j=1}^{+\infty}\overline{\varphi}_j\varphi_j$ as an estimator of $s$
because the series does not converge. This is actually not surprising because we are trying to
estimate infinitely many parameters (the $s_j$) from a finite number of observations. But, for
any finite subset $m$ of $\Bbb{N}^\star$, the estimator $\hat{s}_m=\sum_{j\in m}
\overline{\varphi}_j\varphi_j$ does belong to $\Bbb{L}_2([0,1],dx)$ and
\[
\|\hat{s}_m-s\|^2=\sum_{j\in m}\left(\overline{\varphi}_j-s_j\right)^2+
\sum_{j\not\in m}s_j^2.
\]
If we denote by $|m|$ the cardinality of $m$, we conclude from (\ref{Eq-proj0}) that
\begin{equation}
\Bbb{E}_s\left[\|\hat{s}_m-s\|^2\right]\le2n^{-1}|m|+\|s_m-s\|^2
\quad\mbox{with}\quad s_m=\sum_{j\in m}s_j\varphi_j.
\labe{Eq-proj2}
\end{equation}
Note that $\hat{s}_m$ is not necessarily a genuine estimator, i.e.\ a density, but this is a
minor point since ${\cal S}_2$ is a closed convex subset of $\Bbb{L}_2([0,1],dx)$ on which
we may always project $\hat{s}_m$, getting a genuine estimator which is even closer to $s$
than $\hat{s}_m$.

\subsubsection{Approximate models for nonparametric estimation\labs{D1b}}
The construction of the projection estimator $\hat{s}_m$ can also be interpreted in terms of 
a model since it is actually based on the parametric model
\[ 
S_m=\left\{\left.t=\sum_{j\in m}t_j\varphi_j\,\right|\,t_j\in\Bbb{R}
\mbox{ for }j\in m\right\}.
\]
To build $\hat{s}_m$, we proceed as if $s$ did belong to $S_m$, estimating the $|m|$
unknown parameters $s_j$ for $j\in m$ by their natural estimators $\overline{\varphi}_j$. But
there are three main differences with the classical parametric approach:

i) we do not assume that $s\in S_m$ so that $S_m$ is an {\em approximate model} for the true
density;

ii) apart from some exceptional cases, like histogram estimation, projection estimators are not
maximum likelihood estimators with respect to $S_m$;

iii) there is no asymptotic point of view here and the risk bound (\ref{Eq-proj2}) is valid for
any value of $n$.

The histogram estimator can actually be viewed as a particular projection estimator. With
the notations of Section~\ref{H}, we set $\varphi_j=|I_j|^{-1/2}\1_{I_j}$ for $1\le j\le D$,
we complete this orthonormal family into a basis of $\Bbb{L}_2([0,1],dx)$ and take for $m$
the set $\{1,\dots,D\}$. Then, for $j\in m$,
\[
\overline{\varphi}_j=n^{-1}\sum_{i=1}^n|I_j|^{-1/2}\1_{I_j}(X_i)=n^{-1}|I_j|^{-1/2}N_j
\quad\mbox{and}\quad\sum_{j\in m}\overline{\varphi}_j\varphi_j=\hat{s}_{m}.
\]

\subsection{Approximate models for parametric estimation\labs{D2}}

\subsubsection{Gaussian linear regression\labs{D2a}}
An extremely popular parametric model is {\em Gaussian linear regression}. In this case
we observe $n$ independent variables $X_1,\ldots, X_n$ from the Gaussian linear 
regression set up
\begin{equation}
X_i=\sum_{j=1}^p\beta_jZ^j_i+\sigma\xi_i\quad\mbox{for }1\le i\le n,
\labe{LM1}
\end{equation}
where the random variables $\xi_i$ are i.i.d.\ standard normal while the numbers $Z_i^j$,
$1\le i\le n$ denote the respective deterministic and observable values of some {\em
explanatory variable} $Z^j$. Here, ``variable" is taken in its usual sense of an ``economic
variable" or a ``physical variable".  Practically speaking, $X_i$ corresponds to an observation
in the $i^{\rm{th}}$ experiment and it is assumed that this value depends linearily on the
values $Z_i^j$ of the variables $Z^j, 1\le j\le p$ in this experiment but with some additional
random perturbation represented by the random variable $\sigma\xi_i$. We assume here that
all $p$ parameters $\beta_j$ are unknown but that $\sigma$ is known (this is not usually the
case but will greatly simplify our analysis). This set-up results in a parametric model with $p$
unknown parameters, since the distribution in $\Bbb{R}^n$ of the vector $\bm{X}$ with
coordinates $X_i$ is entirely defined by the parameters $\beta_j$. More precisely, the random
variables $X_1,\ldots,X_n$ are independent with respective normal distributions 
${\cal N}\left(s_i,\sigma^2\right)$ with $s_i=\sum_{j=1}^p\beta_jZ^j_i$. Equivalently
$\bm{X}$ is a Gaussian vector with mean vector $s=(s_i)_{1\le i\le n}$ and covariance matrix
$\sigma^2I_n$ where $I_n$ denotes the identity matrix in $\Bbb{R}^n$. If we denote by
$\bm{Z}^j$ the vector with coordinates $Z^j_i$ and assume that the vectors
$\bm{Z}^j, 1\le j\le p$ span a $p$-dimensional linear space $\overline{S}_p$, which we shall
do, it is equivalent to estimate the parameters $\beta_j$ or the vector $s\in\overline{S}_p$. 

The estimation problem can then be summarized as follows: observing the Gaussian vector
$\bm{X}$ with distribution ${\cal N}\left(s,\sigma^2I_n\right)$ with a known value of
$\sigma$, estimate the parameter $s$ which is assumed to belong to $\overline{S}_p$. This is
a parametric problem similar to those we considered in Section~\ref{C} and it can be solved via
the maximum likelihood method. The density of $\bm{X}$ with respect to the Lebesgue
measure on $\Bbb{R}^n$ and the log-likelihood of $s$ are respectively given by
\[
\frac{1}{\left(2\pi\sigma^2\right)^{n/2}}\exp\left[-\frac{1}{2\sigma^2}\sum_{i=1}^n
(x_i-s_i)^2\right]\quad\mbox{and}\quad
-\frac{n}{2}\log\left(2\pi\sigma^2\right)-\frac{1}{2\sigma^2}\sum_{i=1}^n(X_i-s_i)^2,
\]
so that the maximum likelihood estimator $\hat{s}_p$ over $\overline{S}_p$ is merely the
orthogonal projection of $\bm{X}$ onto $\overline{S}_p$ with risk
$\Bbb{E}_s\left[\|s-\hat{s}_p\|^2\right]=\sigma^2p$. This estimator actually makes sense
even if $s\not\in\overline{S}_p$ since, whatever the true value of $s\in\Bbb{R}^n$,
\begin{equation}
\Bbb{E}_s\left[\|s-\hat{s}_p\|^2\right]=\sigma^2p+\inf_{t\in\overline{S}_p}\|s-t\|^2.
\labe{Eq-Gauss1}
\end{equation}
The risk is the sum of two terms, one which is proportional to the number $p$ of parameters
to be estimated and another one which measures the accuracy of the model $\overline{S}_p$
we use. This second term vanishes when the model is correct (contains $s$). 

\subsubsection{Model choice again\labs{D2b}}
In the classical regression problem, the model $\overline{S}_p$ is assumed to be correct so that
$\Bbb{E}_s\left[\|s-\hat{s}_p\|^2\right]=\sigma^2p$ but this approach leads to two
opposite problems. In order to keep the term $\sigma^2p$ in (\ref{Eq-Gauss1}) small, we
may be tempted to put too few explanatory variables in the model,  omitting some important
ones so that not only $s\not\in\overline{S}_p$ but $\inf_{t\in\overline{S}_p}\|s-t\|^2$ may
be very large, possibly larger than $\sigma^2n$. In this case, it would be wiser to use the largest
possible model $\Bbb{R}^n$ for $s$ and the corresponding m.l.e.\ $\hat{s}=\bm{X}$
resulting in the better risk $\Bbb{E}_s\left[\|s-\hat{s}\|^2\right]=\sigma^2n$. In order to
avoid this difficulty, we may alternatively introduce many explanatory variables $Z^j$ in the
model $\overline{S}_p$. Then even if it is correct, we shall get a large risk bound $\sigma^2p$.
It may then happen that only a small number $q$ of the $p$ explanatory variables
determining the model are really influential. This means that if $\overline{S}_q$ is the linear
span of those $q$ variables, say $Z^1,\ldots,Z^q$, $\inf_{t\in\overline{S}_q}\|s-t\|^2$ is
small. As a consequence, the risk bound of the m.l.e.\ $\hat{s}_q$ with respect to
$\overline{S}_q$, i.e.\ $\sigma^2q+\inf_{t\in\overline{S}_q}\|s-t\|^2$ may be much
smaller than $\sigma^2p$.

These examples show that, even in the parametric case, the use of an approximate model may
be preferable to the use of a correct model, although a grossly wrong model may lead to terrible
results. The choice of a suitable model is therefore crucial: a large model including many
explanatory variables automatically results in a large risk bound due to the component
$\sigma^2p$ of the risk in (\ref{Eq-Gauss1}) while the choice of a too parsimonious model
including only a  limited number of variables may result in a poor estimator based on a grossly
wrong model if we have omitted some very influential variables.

A natural idea to solve this dilemma would be to start with some large family
$\{\overline{S}_m,m\in{\cal M}\}$ of linear models indexed by some set ${\cal M}$ and with
respective dimensions $D_m$. For each of them, the corresponding m.l.e.\ $\hat{s}_m$ (the
projection of $\bm{X}$ onto $\overline{S}_m$) satisfies
\[
\Bbb{E}_s\left[\|s-\hat{s}_m\|^2\right]=\sigma^2D_m+
\inf_{t\in\overline{S}_m}\|s-t\|^2,
\]
and an optimal model $\overline{S}_{\overline{m}}$ is one that minimizes this quantity. But,
as in the case of histograms, this optimal model depends on the unknown parameter $s$ via
$\inf_{t\in\overline{S}_m}\|s-t\|$ so that $\hat{s}_{\overline{m}}$ is an ``oracle", not a
genuine estimator. Since this oracle is not available to the statistician, he has to try an
alternative method and use the observation $\bm{X}$ to build a selection procedure
$\hat{m}(\bm{X})$ of one model $\overline{S}_{\hat{m}}$, estimating $s$ by
$\tilde{s}=\hat{s}_{\hat{m}}$. An ideal {\em model selection procedure} should have the 
performance of an oracle, i.e.\ satisfy
\begin{equation}
\Bbb{E}_s\left[\|s-\tilde{s}\|^2\right]=
\inf_{m\in{\cal M}}\left\{\sigma^2D_m+\inf_{t\in\overline{S}_m}\|s-t\|^2\right\},
\labe{Eq-or1}
\end{equation}
but such a procedure cannot exist and the best that one can expect is to find selection
procedures satisfying a risk bound which is close to (\ref{Eq-or1}).

\Section{Model based statistical estimation\labs{A}}
In three different contexts, namely histogram estimation for densities, projection estimation
for densities and Gaussian linear regression, we have seen that the use of an approximate
model associated with a convenient estimator with values in the model leads to three risk
bounds, namely (\ref{Eq-his6}), (\ref{Eq-proj2}) and (\ref{Eq-Gauss1}), which share the same
structure. These bounds are the sum of two terms, one is the squared distance of the unknown
parameter to the model, the second is proportional to the number of parameters that are
involved in the model. One can therefore wonder to what extent this situation is typical. 

\subsection{A general statistical framework\labs{A0}}
Before we proceed to the solution of the problem, let us make the statistical framework on
which we work somewhat more precise. We observe a random phenomenon
$\bm{X}(\omega)$ (real variable, vector, sequence, process, set, \ldots) from the abstract
probability space $(\Omega,{\cal A},\Bbb{P})$ with values in the measurable set $(\Xi,{\cal
X})$ and with unknown probability distribution
$P_{\bm{X}}$ on $(\Xi,{\cal X})$ given by 
\[ 
P_{\bm{X}}[A]=\Bbb{P}\left[\bm{X}^{-1}(A)\right]=\Bbb{P}\left[\{\omega\in\Omega\,|
\,\bm{X}(\omega)\in A\}\right]\quad\mbox{for all }A\in{\cal X}.
\] 
The purpose of statistical estimation is to get some information on this distribution from one
observation $\bm{X}(\omega)$ of the phenomenon. We assume that $P_{\bm{X}}$ belongs
to some given subset ${\cal P}=\{P_t, t\in M\}$ of the set of all distributions on $(\Xi,{\cal X})$, where $M$ denotes a one-to-one parametrization of ${\cal P}$. We moreover assume
that $M$ is a metric space with a distance $d$. Therefore $P_{\bm{X}}=P_s$ for some $s\in
M$ and we want to {\em estimate} $P_s$, or equivalently $s$, in view of this one-to-one
correspondence which also allows us to consider $d$ as a distance on ${\cal P}$ as well.
As in Section~\ref{H2b}, we look for an estimator of $s$, i.e.\ a measurable mapping
$\hat{s}$ from $(\Xi,{\cal X})$ to $M$ (with its Borel $\sigma$-algebra) such that
$\hat{s}(\bm{X})$ provides a good approximation of the unknown value $s$. Such a
mapping is called an {\em estimator} of $s$. We measure the performance of the estimator
$\hat{s}(\bm{X})$ via its {\em quadratic risk}
\begin{equation}
R\left(\hat{s},s\right)=\Bbb{E}_s\left[d^2\left(s,\hat{s}\right)\right].
\labe{Eq-risk1}
\end{equation}
There is a very large number of possibilities for the choice of ${\cal P}$ depending on the
structure of $(\Xi,{\cal X})$ and the problem we have to solve. In this paper we
focus on the two particular but typical examples that we considered earlier, namely the
density estimation problem and the Gaussian regression problem which amounts to the
estimation of the mean of a Gaussian vector. In both cases $\Xi=E^n$ is a product space with a
product $\sigma$-algebra ${\cal X}={\cal E}^{\otimes n}$ so that $\bm{X}$ is the vector
$(X_1,\ldots,X_n)$ and the $X_i$ are random variables with values in $(E,{\cal E})$. 

\paragraph{Density estimation}
For the {\em density estimation problem} we are given some reference measure $\nu$ on
$(E,{\cal E})$ and we assume that the $X_i$ are  i.i.d.\ random variables with a density $s$
with respect to $\nu$, in which case $M$ can be chosen as the set of all densities with
respect to $\nu$, i.e.\ the subset of $\Bbb{L}_1(\nu)$ of nonnegative functions which
integrate to one. Such a situation occurs when one replicates the same experiment $n$ times
under identical conditions and assumes that each experiment has no influence on the others,
for instance when we observe the successive outcomes of a ``roulette" game. Then, for each
$t\in M$, $P_t$ has the density $\prod_{i=1}^n t(x_i)$ with respect to
$\mu=\nu^{\otimes n}$.

\paragraph{Gaussian regression}
This is the case that we considered in Section~\ref{D2} with $(E,{\cal E})$ being the real
line with its Borel $\sigma$-algebra. Here $\bm{X}$ is a Gaussian vector in $\Bbb{R}^n$
with known covariance matrix $\sigma^2 I_n$. Then $M=\Bbb{R}^n$ and
$t=(t_1,\ldots,t_n)\in M$ is the unknown mean vector of the Gaussian distribution $P_t={\cal
N}\left(t,\sigma^2 I_n\right)$ with density
\begin{equation}
g_t(x)=\frac{1}{\left(2\pi\sigma^2\right)^{n/2}}\exp\left[-\frac{1}{2\sigma^2}
\sum_{i=1}^n(x_i-t_i)^2\right],
\labe{Eq-Gauss0}
\end{equation}
with respect to the Lebesgue measure $\mu$ on $\Bbb{R}^n$.

\subsection{Two point parameter sets\labs{A1}}
Before we come to the general situation, it will be useful to analyze a special, quite irrealistic,
but very simple case. Let us make the extra assumption that $s$ belongs to the smallest
possible parameter set, i.e.\ a subset $S$ of $M$ containing only two elements $v$ and
$u$. Note that the statistical problem would be void if $S$ contained only one point
since $s$ would then be known.

A solution to this estimation problem is provided by the maximum likelihood method
described in Section~\ref{C2}. Let $\mu$ be any measure dominating both $P_v$ and $P_u$
($P_v+P_u$ would do) and denote by $g_v$ and $g_u$ the respective densities of $P_v$ and
$P_u$ with respect to $\mu$. Then define an estimator $\hat{\varphi}(\bm{X})$ with values in
$S$ by 
\begin{equation}
\hat{\varphi}(\bm{X})=\left\{\begin{array}{lll}v&\mbox{if}&g_v(\bm{X})>
g_u(\bm{X});\\u&\mbox{if}&g_u(\bm{X})>g_v(\bm{X}). \end{array}\right.
\labe{Eq-mle1}
\end{equation}
Take any decision you like in case of equality. If $s=v$, we get
\[
R\left(\hat{\varphi},s\right)=d^2(v,u)\Bbb{P}_v[\hat{\varphi}=u]\le
d^2(v,u)\Bbb{P}_v[g_u(\bm{X})\ge g_v(\bm{X})].
\]
Since the distribution of $\bm{X}$ is $P_s=P_v=g_v\cdot\mu$, $\Bbb{P}_v[g_v(\bm{X})>0]=1$
and
\begin{eqnarray*}
\Bbb{P}_v[g_u(\bm{X})\ge
g_v(\bm{X})]&=&\Bbb{P}_v\left[\sqrt{g_u(\bm{X})/g_v(\bm{X})}\ge1\right]\;\;\le\;\;
\Bbb{E}_v\left[\sqrt{g_u(\bm{X})/g_v(\bm{X})}\right]\\&=&\int_\Xi
\sqrt{g_u(x)/g_v(x)}\,g_v(x)\,d\mu(x)\;\;=\;\;\int_\Xi\sqrt{g_u(x)g_v(x)}\,d\mu(x).
\end{eqnarray*}
Hence $R\left(\hat{\varphi},s\right)\le d^2(v,u)\rho(P_v,P_u)$ with
\begin{equation}
\rho(P_v,P_u)=\rho(P_u,P_v)=
\int_\Xi\sqrt{\frac{dP_u}{d\mu}(x)\frac{dP_v}{d\mu}(x)}\,d\mu(x).
\labe{Eq-rho0}
\end{equation}
It is easily seen that the definition of $\rho(P_v,P_u)$ via (\ref{Eq-rho0}) is independent of the
choice of the dominating measure $\mu$. Since the same risk bound holds when $s=u$, we
finally get
\begin{equation}
\sup_{s\in\{u,v\}}R\left(\hat{\varphi},s\right)\le d^2(v,u)\rho(P_v,P_u).
\labe{Eq-rb2}
\end{equation}
This bound demonstrates the importance of the so-called {\em Hellinger affinity} $\rho(P,Q)$
between two probability measures $P$ and $Q$. It satisfies in particular by the
Cauchy-Schwarz Inequality and the Fubini Theorem
\begin{equation}
0\le\rho(P,Q)\le 1\qquad\mbox{and}\qquad
\rho\left(P^{\otimes n},Q^{\otimes n}\right)=\rho^n(P,Q).
\labe{Eq-rho}
\end{equation}
It is, moreover, closely related to a well-known distance between probabilities, the 
{\em Hellinger distance} $h$ defined by 
\begin{equation}
h^2(P,Q)=\frac{1}{2}\int_\Xi\left(\sqrt{\frac{dP}{d\mu}(x)}-
\sqrt{\frac{dQ}{d\mu}(x)}\right)^2d\mu(x)=1-\rho(P,Q).
\labe{Eq-hell}
\end{equation}
The Hellinger distance is merely the $\Bbb{L}_2(\mu)$-distance between the square roots of the
densities with respect to any dominating measure $\mu$ (and actually  independent of $\mu$).
Here, we follow Le Cam who normalizes the integral so that the Hellinger distance has range
$[0,1]$. An alternative definition is without the factor $1/2$ in (\ref{Eq-hell}). He also showed
in Le Cam (1973) that
\begin{equation}
\rho(P,Q)\ge\int_\Xi\inf\left\{\frac{dP}{d\mu}(x);\frac{dQ}{d\mu}(x)\right\}d\mu(x)
\ge1-\sqrt{1-\rho^2(P,Q)}.
\labe{Eq-ropi}
\end{equation}

It is easy to compute $\rho(P_v,P_u)$ for our two special frameworks. In the case of Gaussian
distributions $P_u={\cal N}\left(u,\sigma^2 I_n\right)$ and
$P_v={\cal N}\left(v,\sigma^2 I_n\right)$, we get
\[
\rho(P_v,P_u)=\exp\left[-\|v-u\|^2/\left(8\sigma^2\right)\right],
\]
so that $[-\log\rho]^{1/2}$ is a multiple of the Euclidian distance between parameters,
modulo the identification of $t$ and $P_t$. Note that, in general, $[-\log\rho]^{1/2}$ is not a
distance since it may be infinite and does not satisfy the triangle inequality. Setting
$d(v,u)=\|v-u\|$, (\ref{Eq-rb2}) becomes
\[
\sup_{s\in\{u,v\}}R\left(\hat{\varphi},s\right)\le
\|v-u\|^2\exp\left[-\|v-u\|^2/\left(8\sigma^2\right)\right]\le8e^{-1}\sigma^2,
\]
independently of $v$ and $u$. In the i.i.d.\ case, we use the Hellinger distance to define the
risk, setting $d(u,v)=h(u,v)=h(P_u,P_v)$ and (\ref{Eq-rb2}) becomes, whatever the densities $v$
and
$u$,
\[
\sup_{s\in\{u,v\}}R\left(\hat{\varphi},s\right)\le h^2(v,u)\left[1-h^2(v,u)\right]^n\le
n^n(n+1)^{-(n+1)}\le (ne)^{-1}.
\]

\subsection{Two point models for the Gaussian framework\labs{A2}}
As we pointed out at the beginning of the last section, assuming that $s$ is either $v$ or $u$
is definitely irrealistic. A more realistic problem would rather be as follows: $s$ is unknown
but we believe that one of two different situations can occur implying that $s$ is close (not
necessarily equal) to either $v$ or $u$. Then it seems natural to use $S=\{v,u\}$ as an
approximate model for $s$ and just proceed as before, using the estimator
$\hat{\varphi}(\bm{X})$ defined by (\ref{Eq-mle1}). We can then try to mimic the proof
which lead to (\ref{Eq-rb2}), apart from the fact that the argument leading to
\[
\Bbb{P}_v[g_u(\bm{X})\ge g_v(\bm{X})]\le\rho(P_v,P_u)=
\exp\left[-\|v-u\|^2/\left(8\sigma^2\right)\right]
\]
then fails. One can instead prove the following result (Birg\'e, 2006).
%
\begin{proposition}\lab{P-Gtest}
Let $P_t$ denote the Gaussian distribution ${\cal N}\left(t,\sigma^2I_n\right)$ in
$\Bbb{R}^n$. If $\bm{X}$ is a Gaussian vector with distribution $P_s$ and
$\|s-v\|\le\|v-u\|/6$, then
\[
\Bbb{P}_s[g_u(\bm{X})\ge g_v(\bm{X})]
\le\exp\left[-\|v-u\|^2/\left(24\sigma^2\right)\right].
\]
\end{proposition}
We can then proceed as before and conclude that, if $\|s-v\|\le\|v-u\|/6$, then
\begin{eqnarray*}
R\left(\hat{\varphi},s\right)&\le&2\left(\|s-v\|^2+
\Bbb{E}_s\left[\|\hat{\varphi}-v\|^2\right]\right)\\&\le
&2\left(\|s-v\|^2+\|v-u\|^2\exp\left[-\|v-u\|^2/\left(24\sigma^2\right)\right]\right)
\\&\le&2\|s-v\|^2+48e^{-1}\sigma^2.
\end{eqnarray*}
A similar bound holds with $u$ replacing $v$ if $\|s-u\|\le\|v-u\|/6$. Finally, if   
$\min\{\|s-v\|,\|s-u\|\}>\|v-u\|/6$, since $\hat{\varphi}$ is either $v$ or $u$,
\begin{eqnarray*}
R\left(\hat{\varphi},s\right)&\le&(\max\{\|s-v\|,\|s-u\|\})^2\\&\le&
(\min\{\|s-v\|,\|s-u\|\}+\|v-u\|)^2\\&\le&49(\min\{\|s-v\|,\|s-u\|\})^2.
\end{eqnarray*}
We finally conclude that, whatever $s\in M$, even if our initial assumption that $s$ is close
to $S$ is wrong,
\[
R\left(\hat{\varphi},s\right)\le 48e^{-1}\sigma^2+49\inf_{t\in S}\|s-t\|^2,
\]
which, apart from the constants, is similar to (\ref{Eq-Gauss1}).

\subsection{General models for the Gaussian  framework\labs{A3}}

\subsubsection{Linear models\labs{A3a}}
Instead of assuming that $s$ is close to a two-points set, let us now assume that it is close to
some $D$-dimensional linear subspace $V$ of $\Bbb{R}^n$ ($D>0$). Choose some
$\lambda\ge4\sqrt{3}\sigma$ and, identifying $V$ to $\Bbb{R}^{D}$ via some
orthonormal basis, consider the lattice $S=(2\lambda\Bbb{Z})^{D}\subset V$. The maximum
likelihood estimator $\hat{s}(\bm{X})$ with respect to $S$ is given by $\hat{s}(\bm{X})=
\argmax_{t\in S}g_t(\bm{X})$. Its unicity follows from the facts that $S$ is countable and
$\Bbb{P}_s[g_t(\bm{X})=g_u\bm{X})]=0$ for each pair $(t,u)\in S^2$ such that $t\ne u$. As to
its existence (with probability one), it is a consequence of the following result.
%
\begin{proposition}\lab{P-mod1}
For $s$ an arbitrary point in $M=\Bbb{R}^n$, $s'\in S$, and 
\begin{equation}
y\ge y_0=\max\left\{\lambda\sqrt{2D},6\|s'-s\|\right\},
\labe{EqM3}
\end{equation}
then
\begin{equation}
\Bbb{P}_s\left[\exists t\in S\mbox{ with }\|s'-t\|\ge y\mbox{ and }
g_t(\bm{X})\ge g_{s'}(\bm{X})\right]\le1.14\exp\left[-\frac{y^2}{48\sigma^2}\right].
\labe{Eq-expbou}
\end{equation}
\end{proposition}
%
\noindent{\em Proof:}
Let $S_k=\left\{t\in S\,\left|\,2^{k/2}y\le\|s'-t\|<2^{(k+1)/2}y\right.\right\}$ with
cardinality  $|S_k|$. If we denote by $P(y)$ the left-hand side of (\ref{Eq-expbou}), we get
\begin{equation}
P(y)\le\sum_{k=0}^{+\infty}\Bbb{P}_s[\exists t\in S_k\mbox{ with }g_t(\bm{X})\ge
g_{s'}(\bm{X})]\nonumber\le\sum_{k=0}^{+\infty}|S_k|\sup_{t\in S_k}
\Bbb{P}_s[g_t(\bm{X})\ge g_{s'}(\bm{X})].
\labe{EqM2}
\end{equation}
Since, for $t\in S_k$, $\|s'-t\|\ge2^{k/2}y\ge6\|s'-s\|$, we may apply
Proposition~\ref{P-Gtest} to get 
\begin{equation}
\sup_{t\in S_k}\Bbb{P}_s[g_t(\bm{X})\ge g_{s'}(\bm{X})]\le
\exp\left[-2^ky^2/\left(24\sigma^2\right)\right].
\labe{EqM1}
\end{equation}
Moreover, for any ball ${\cal B}(s',r)$ with center $s'$ and radius $r=x\lambda\sqrt{D}$ with
$x\ge2$,
\begin{equation}
|S\cap{\cal B}(s',r)|<\exp\left[x^2D/2\right].
\labe{Eq-app2}
\end{equation}
To prove this, we apply the next inequality which follows from a comparison of the volumes of
cubes and balls in $\Bbb{R}^D$ as in the proof of Lemma~2 from Birg\'e and Massart (1998).
\[
|S\cap{\cal B}(s',r)|\le\frac{(\pi e/2)^{D/2}}{\sqrt{\pi D}}
\left(\frac{r}{\lambda\sqrt{D}}+1\right)^{D}<\exp[D(0.73+\log(x+1))].
\]
We then get (\ref{Eq-app2}) since $x\ge2$. Applying it with $r=2^{(k+1)/2}y\ge2^{1+k/2}
\lambda\sqrt{D}$ by (\ref{EqM3}), leads to $|S_k|\le\exp\left[2^k(y/\lambda)^2\right]$.
Together with (\ref{EqM2}) and (\ref{EqM1}), this shows that
\begin{eqnarray*}
P(y)&\le&\sum_{k=0}^{+\infty}\exp\left[2^k\frac{y^2}{\lambda^2}
-2^k\frac{y^2}{24\sigma^2}\right]\;\;\le\;\;\sum_{k=0}^{+\infty}
\exp\left[-2^k\frac{y^2}{48\sigma^2}\right]\\&=&
\exp\left[-\frac{y^2}{48\sigma^2}\right]
\sum_{k=0}^{+\infty}\exp\left[-\frac{y^2}{48\sigma^2}\left(2^k-1\right)\right].
\end{eqnarray*}
The conclusion follows from the fact that
$y^2\ge y_0^2\ge2\lambda^2D\ge2\lambda^2\ge96\sigma^2$.\cqfd\vspace{2mm}\\
Proposition~\ref{P-mod1} implies that, for $y\ge y_0$, there exists a set $\Omega_y\subset
\Omega$ with
$\Bbb{P}_s(\Omega_y)\ge1-1.14\exp\left[-y^2/\left(48\sigma^2\right)\right]$ and such
that, for $\omega\in\Omega_y$, the function $t\mapsto g_t(\bm{X}(\omega))$ has a
maximum in the ball ${\cal B}(s',y)$. This shows that, if $\omega\in\Omega_y$, the m.l.e.\
$\hat{s}(\bm{X}(\omega))$ exists and satisfies $\|\hat{s}(\bm{X})-s'\|\le y$. As a
consequence, the m.l.e.\ $\hat{s}(\bm{X})$ exists a.s.\ and
\begin{eqnarray*}
\Bbb{E}_s\left[\|\hat{s}(\bm{X})-s'\|^2\right]&=&\int_0^{+\infty}\Bbb{P}_s
\left[\|\hat{s}(\bm{X})-s'\|^2\ge z\right]dz\\&\le&y_0^2+\int_{y_0^2}^{+\infty}
\Bbb{P}_s\left[\|\hat{s}(\bm{X})-s'\|\ge\sqrt{z}\right]dz\\&\le&y_0^2+1.14
\int_{y_0^2}^{+\infty}\exp\left[-\frac{z}{48\sigma^2}\right]dz\\&=&y_0^2+1.14\times
48\sigma^2\exp\left[-y_0^2/\left(48\sigma^2\right)\right]\\&\le&y_0^2+
55e^{-2}\sigma^2.
\end{eqnarray*}
Then
\begin{eqnarray*}
\Bbb{E}_s\left[\|\hat{s}(\bm{X})-s\|^2\right]&\le&2\left[\|s-s'\|^2+
\Bbb{E}_s\left[\|\hat{s}(\bm{X})-s'\|^2\right]\right]\\&\le&2\left[\|s-s'\|^2+y_0^2+55
e^{-2}\sigma^2\right]\\&\le&2\left[37\|s-s'\|^2+2\lambda^2D+55e^{-2}\sigma^2\right].
\end{eqnarray*}
Note that the construction of $S$ as a lattice in $V$ implies that any point in $V$ is at a
distance of some point in $S$ not larger than $\lambda\sqrt{D}$ which means that one can
choose $s'$ in such a way that $\|s-s'\|\le\inf_{t\in V}\|s-t\|+\lambda\sqrt{D}$. With
such a choice for $s'$, we get
\[
\Bbb{E}_s\left[\|\hat{s}(\bm{X})-s\|^2\right]\le2\left[74\inf_{t\in V}\|s-t\|^2+
76\lambda^2D+55e^{-2}\sigma^2\right].
\]
Setting $\lambda$ to its minimum value $4\sqrt{3}\sigma$, we conclude, since $D\ge1$,
that
\begin{equation}
\Bbb{E}_s\left[\|\hat{s}(\bm{X})-s\|^2\right]\le148\inf_{t\in V}\|s-t\|^2+
7311\sigma^2D.
\labe{Eq-M6}
\end{equation}

\subsubsection{General models with finite metric dimension\labs{A3b}}
Note that, apart from the huge constants that we actually did not try to optimize in order to
keep the computations as simple as possible, (\ref{Eq-M6}) is quite similar to (\ref{Eq-Gauss1}),
although we actually used a different estimation procedure, and also a different method of
proof which has an important advantage: it did not make any use of the fact that $V$ is a linear
space. What we actually used are the metric properties of the $D$-dimensional linear subspace
$V$ of $M=\Bbb{R}^n$, which can be summarized as follows.
%
\begin{property}
Whatever $\eta>0$, one can find a subset $S$ of $M$ such that:

i) for each $t\in V$ there exists some $t'\in S$ with $\|t-t'\|\le\eta$;

ii) for any ball ${\cal B}(t,x\eta)$ with center $t\in M$ and radius $x\eta$,
\[
|S\cap{\cal B}(t,x\eta)|\le\exp\left[x^2D/2\right]\quad\mbox{for }x\ge2.
\]
\end{property}
In the previous example we simply defined $S$ so that
$\eta=\lambda\sqrt{D}=4\sigma\sqrt{3D}$.

The fact that the previous property of $V$ was a key argument in the proof motivates the
following general definition.

\begin{definition}\lab{D-mdim}
Let $\overline{S}$ be a subset of some metric space $(M,d)$ and $\overline{D}$ be some real
number $\ge1/2$. We say that $\overline{S}$ has {\em a finite metric dimension bounded by
$\overline{D}$} if, for every $\eta>0$, one can find a subset $S_\eta$ of $M$ such that:

i) for each $t\in\overline{S}$ there exists some $t'\in S_\eta$ with $d(t,t')\le\eta$ (we say
that $S_\eta$ is {\em an $\eta$-net for $\overline{S}$});

ii) for any ball ${\cal B}(t,x\eta)$ with center $t\in M$ and radius $x\eta$,
\[
|S_\eta\cap{\cal B}(t,x\eta)|\le\exp\left[x^2\overline{D}\right]\quad\mbox{for }x\ge2.
\]
\end{definition}
Note that any subset of $\overline{S}$ also has a finite metric dimension bounded by
$\overline{D}$. It follows from the Property~P that a $D$-dimensional linear subspace
of a Euclidean space has a metric dimension bounded by $D/2$. Note that, apart from the
factor 1/2, this result cannot be improved in view of the following lower bound for the metric
dimension of a $D$-dimensional ball.
%
\begin{lemma}\lab{L-infdim}
Let $\overline{S}$ be a ball of the metric space $(M,d)$ which is isometric to a ball in the
Euclidean space $\Bbb{R}^D$. Then a bound $\overline{D}$ for its metric dimension cannot be
smaller than $D/13$.
\end{lemma}
%
\noindent{\em Proof:}
Let $\overline{S}={\cal B}(t,r)$ have a finite metric dimension bounded by $\overline{D}$
and $\eta<r/3$. One can find $S_\eta$ in $M$ which is an $\eta$-net for $\overline{S}$
and such that $N=|S_\eta\cap{\cal B}(t,3\eta)|\le\exp\left[9\overline{D}\right]$.
Moreover, $S_\eta$ is also an $\eta$-net for ${\cal B}(t,2\eta)$ so that ${\cal B}(t,2\eta)$ can
be covered by the $N$ balls with radius $\eta$ and centers in $S_\eta\cap{\cal B}(t,3\eta)$.
Since ${\cal B}(t,3\eta)\subset\overline{S}$ we can use the isometry to show, comparing the
volumes of the balls, that $N\ge2^D$ so that $9\overline{D}\ge D\log2$ and the conclusion
follows.\cqfd\vspace{2mm}

Introducing Definition~\ref{D-mdim} in the proof of Proposition~\ref{P-mod1}, we get the
following result.
%
\begin{theorem}\lab{T-Gmodel}
Let $\bm{X}$ be a Gaussian vector in $\Bbb{R}^n$ with unknown mean $s$ and known
covariance matrix $\sigma^2I_n$. Let $\overline{S}$ be a subset of the Euclidean space
$\Bbb{R}^n$ with a finite metric dimension bounded by $\overline{D}$. Then one can
build an estimator $\hat{s}_{\overline{S}}(\bm{X})$ of $s$ such that, for some universal
constant $C$ (independent of $s$, $n$ and $\overline{S}$),
\begin{equation}
\Bbb{E}_s\left[\|\hat{s}_{\overline{S}}(\bm{X})-s\|^2\right]\le
C\left[\inf_{t\in\overline{S}}\|s-t\|^2+\sigma^2\overline{D}\right].
\labe{Eq-rbg}
\end{equation}
\end{theorem}
This theorem implies that we can use for models non-linear sets that have a finite metric
dimension. In particular, various types of manifolds could be used as models. To build the
estimator $\hat{s}_{\overline{S}}(\bm{X})$, we set $\eta=4\sigma\sqrt{6\overline{D}}$ and
choose an $\eta$-net $S_\eta$ for $\overline{S}$ satisfying the properties of
Definition~\ref{D-mdim}. Then we take for $\hat{s}_{\overline{S}}(\bm{X})$ the m.l.e.\ with
respect to $S_\eta$.

\subsection{Density estimation\labs{A4}}
When we want to extend the results obtained for the Gaussian framework to density
estimation we encounter new difficulties. The two key arguments used in the proof of
Proposition~\ref{P-mod1} are that $V$ has a finite metric dimension and
Proposition~\ref{P-Gtest}. For i.i.d.\ observations $X_1,\ldots,X_n$ with density $s$ and in
view of the fact that
\[
\Bbb{P}_s\left[\prod_{i=1}^nu(X_i)\ge\prod_{i=1}^nv(X_i)\right]\le
\exp\left[-nh^2(u,v)\right]\quad\mbox{if }s=v,
\]
an analogous result would be as follows:
%
\begin{conjecture}\lab{Z-ftest}
Let $X_1,\ldots,X_n$ be i.i.d.\ random variables with an unknown density $s$ with respect
to some measure  $\nu$ on $(E,{\cal E})$. There exist two constants $\kappa\ge2$ and
$A>0$ such that, whatever the densities $u,v$  on $(E,{\cal E})$ such that
$h(s,v)\le\kappa^{-1}h(u,v)$, then
\[
\Bbb{P}_s\left[\prod_{i=1}^nu(X_i)\ge\prod_{i=1}^nv(X_i)\right]
\le\exp\left[-Anh^2(u,v)\right].
\]
\end{conjecture}
If this conjecture were true one could mimic the proof for the Gaussian case, starting from a
subset $\overline{S}$ of the metric space $(M,h)$ with finite metric dimension, choosing a
suitable $\eta$-net $S_\eta$ for $\overline{S}$ and computing the m.l.e.\ with respect to
$S_\eta$ to get an analogue of Theorem~\ref{T-Gmodel}. Unfortunately Conjecture~C is
wrong and, as a consequence, one can find stuations in the i.i.d.\ framework where the
m.l.e.\ with respect to $S_\eta$ does not behave at all as expected. To get an analogue of
Theorem~\ref{T-Gmodel} for density estimation, one cannot work with the maximum
likelihood method any more. An alternative method that allows to deal with the problem of
density estimation has been proposed by Le Cam (1973 and 1975) who also introduced a notion
of metric dimension, and then extended by the present author in Birg\'e (1983 and 1984). In the
sequel, we shall follow the generalized approach of Birg\'e (2006) from which we borrow this
substitute to Conjecture~C:
%
\begin{proposition}\lab{P-vtest}
Let $X_1,\ldots,X_n$ be i.i.d.\ random variables with an unknown density $s$ with respect
to some measure  $\nu$ on $(E,{\cal E})$. Whatever the densities $u,v$, one can design a
procedure $\varphi_{u,v}(X_1,\ldots,X_n)$ with values in $\{u,v\}$ and such that 
\[
\Bbb{P}_s\left[\varphi_{u,v}(X_1,\ldots,X_n)=u\right]
\le\exp\left[-(n/4)h^2(u,v)\right]\quad\mbox{if }h(s,v)\le h(u,v)/4;
\]
\[
\Bbb{P}_s\left[\varphi_{u,v}(X_1,\ldots,X_n)=v\right]
\le\exp\left[-(n/4)h^2(u,v)\right]\quad\mbox{if }h(s,u)\le h(u,v)/4.
\]
\end{proposition}
The main difference with Conjecture~C lies in the fact that the procedure $\varphi_{u,v}$
does not choose between $u$ and $v$ by merely comparing $\prod_{i=1}^nu(X_i)$ and
$\prod_{i=1}^nv(X_i)$. It is more complicated. This implies that, in this case, we have to
design a new estimator $\hat{s}_{\overline{S}}(\bm{X})$, based on
Proposition~\ref{P-vtest}, to replace the m.l.e.. The construction of this estimator is more
complicated than that of the m.l.e.\ and we shall not describe it here. The following analogue
of Theorem~\ref{T-Gmodel} is proved in Birg\'e (2006).
%
\begin{theorem}\lab{T-Dmodel}
Let $\bm{X}=(X_1,\ldots,X_n)$ be an i.i.d.\ sample with unknown density $s$ with respect
to some measure  $\nu$ on $(E,{\cal E})$ and $(M,h)$ be the metric space of all such
densities with Hellinger distance. Let $\overline{S}$ be a subset of $(M,h)$ with a finite
metric dimension bounded by $\overline{D}$. Then one can build an estimator
$\hat{s}_{\overline{S}}(X_1,\ldots,X_n)$ of $s$ such that, for some universal constant $C$,
\begin{equation}
\Bbb{E}_s\left[h^2\left(\hat{s}_{\overline{S}},s\right)\right]\le
C\left[\inf_{t\in\overline{S}}h^2(s,t)+n^{-1}\overline{D}\right].
\labe{Eq-rbd}
\end{equation}
\end{theorem}
Analogues of Proposition~\ref{P-vtest} do hold for various statistical frameworks, although
not all. Additional examples are to be found in Birg\'e (2004 and 2006). For each such case,
one can, starting from a model $\overline{S}$ with finite metric dimension bounded by
$\overline{D}$, design a suitable estimator $\hat{s}_{\overline{S}}(\bm{X})$ and then get an
analogue of Theorem~\ref{T-Dmodel}. Within the general framework of Section~\ref{A0}, the
resulting risk bound takes the following form:
\begin{equation}
\Bbb{E}_s\left[d^2\left(\hat{s}_{\overline{S}},s\right)\right]\le
C_1\inf_{t\in\overline{S}}d^2(s,t)+C_2\overline{D}\quad\mbox{for all }s\in M,
\labe{Eq-gen1}
\end{equation}
where the constants $C_1$ and $C_2$ depend on the corresponding statistical framework ---
compare with (\ref{Eq-rbg}) and (\ref{Eq-rbd}) --- but not on $s$ or $\overline{S}$. The main
task is indeed to prove the proper alternative to Proposition~\ref{P-vtest}. Once this has been
done, (\ref{Eq-gen1}) follows more or less straightforwardly.

To what extent can maximum likelihood or related estimators provide bounds of the form
(\ref{Eq-gen1}) has been studied in various papers among which van de Geer (1990, 1993,
1995 and 2000), Shen and Wong (1994) and Wong and Shen (1995), Birg\'e and Massart
(1993 and 1998), Gy\"orfi, Kohler, Kry\.zak and Walk (2002) and Massart (2006).

\Section{Model selection\labs{S}}
Let us consider a statistical framework for which an analogue of Proposition~\ref{P-vtest}
holds so that any model $\overline{S}$ with finite metric dimension bounded by
$\overline{D}$ provides an estimator $\hat{s}_{\overline{S}}(\bm{X})$ with a risk
bounded by (\ref{Eq-gen1}). Then the quality of a given model $\overline{S}$ for estimating
$s$ can be measured by the right-hand side of (\ref{Eq-gen1}). Since this quality depends on
the unknown $s$ via the approximation term $\inf_{t\in\overline{S}}d^2(s,t)$, we cannot
know it. Introducing a large family $\{\overline{S}_m, m\in{\cal M}\}$ of models, each one 
with finite metric  dimension bounded by $\overline{D}_m$, instead of one single model, gives
more chance to get an estimator $\hat{s}_m=\hat{s}_{\overline{S}_m}$ in the family with the
smaller risk bound $\inf_{m\in{\cal M}}\left\{C_1\inf_{t\in\overline{S}}d^2(s,t)+
C_2\overline{D}_m\right\}$. Since we do not know which estimator reaches this bound, the
challenge of {\em model selection} is to design a random choice $\hat{m}(\bm{X})$ of $m$ such
that the corresponding estimator $\hat{s}_{\hat{m}}$ approximately reaches this optimal risk,
i.e.\ satisfies
\begin{equation}
\Bbb{E}_s\left[d^2(s,\hat{s}_{\hat{m}})\right]\le C
\inf_{m\in{\cal M}}\left\{C_1\inf_{t\in\overline{S}_m}d^2(s,t)+C_2\overline{D}_m\right\},
\labe{Eq-topt}
\end{equation}
for some constant $C$ independent of $s$ and the family of models.

\subsection{Some natural limitations to the performances of model selection\labs{S1}}
Let us show here, in the context of Gaussian regression, that getting a bound like
(\ref{Eq-topt}) for arbitrary families of models is definitely too optimistic. If, in this context,
(\ref{Eq-topt}) were true, we would be able to design a model selection procedure $\hat{m}$
satisfying, in view of (\ref{Eq-M6})
\begin{equation}
\Bbb{E}_s\left[\|s-\hat{s}_{\hat{m}}\|^2\right]\le C'\inf_{m\in{\cal M}}
\left\{\sigma^2\overline{D}_m+\inf_{t\in\overline{S}_m}\|s-t\|^2\right\},
\labe{Eq-wrb}
\end{equation}
for some universal constant $C'$, independent of $s$, $n$ and the family of models. It is
not difficult to see that this is impossible, even if we restrict ourselves to countable families of
models. Indeed, if (\ref{Eq-wrb}) were true, we could choose for $\{\overline{S}_m,
m\in{\cal M}\}$ a countable family of one-dimensional linear spaces such that each point
$s\in\Bbb{R}^n$ could be approximated by one space in the family with arbitrary accuracy.
We would then get $\overline{D}_m=1/2$ for each $m$ and (\ref{Eq-wrb}) would imply that
\[
\Bbb{E}_s\left[\|s-\hat{s}_{\hat{m}}\|^2\right]\le C'\sigma^2/2
\quad\mbox{for all }s\in\Bbb{R}^n.
\]
But it is known that the best bound one can expect for any estimator $\hat{s}$ uniformly
with respect to $s\in\Bbb{R}^n$ is
\[
\sup_{s\in\Bbb{R}^n}\Bbb{E}_s\left[\|s-\hat{s}\|^2\right]=n\sigma^2,
\]
which contradicts the fact that $C'$ should be a universal constant. One actually has to pay a
price for using many models simultaneously and, as we shall see, this price depends on the
{\em complexity }(with a suitable sense) of the chosen family of models.

\subsection{Risk bounds for model selection\labs{S2}}

\subsubsection{The main theorems\labs{S2a}}
We shall not get here into the details of the construction of the selection procedure that we use
but content ourselves to give the main results and analyze their consequences. A key idea for
the construction appeared in Barron and Cover (1991). Further approaches to selection
procedures have been developed in Barron, Birg\'e and Massart (1999), Birg\'e and
Massart (1997 and 2001), van de Geer (2000), Gy\"orfi, Kohler, Kry\.zak and Walk (2002) and
Massart (2006) who provides an extensive list of references. We follow here the approach based
on dimension from Birg\'e (2006), providing hereafter two theorems corresponding to our two
problems of interest, Gaussian regression and density estimation. In both cases, the construction
of the estimators requires the introduction of a family of positive weights
$\{\Delta_m, m\in{\cal M}\}$, to be chosen by the statistician and satisfying the condition
\begin{equation}
\sum_{m\in{\cal M}}\exp\left[-\Delta_m\right]\le1.
\labe{Eq-S}
\end{equation}
In case of equality in (\ref{Eq-S}), the family $\{q_m\}_{m\in{\cal M}}$ with
$q_m=\exp\left[-\Delta_m\right]$ defines a probability $Q$ on the family of models and
choosing a large value for $\Delta_m$ means putting a small probability on the model
$\overline{S}_m$. One can then see $q_m$ as a probability that the statistician puts on
$\overline{S}_m$ and which influences the result of the estimation procedure, as shown by
the next theorems. Such an interpretation of the weights $\Delta_m$ corresponds to the
so-called {\em Bayesian} point of view. A detailed analysis of this interpretation can be
found in Birg\'e and Massart (2001, Sect.~3.4).
%
\begin{theorem}\lab{T-GP}
Let $\bm{X}$ be a Gaussian vector in $\Bbb{R}^n$ with unknown mean $s$ and known
covariance matrix $\sigma^2I_n$. Let $\{\overline{S}_m, m\in{\cal M}\}$ be a finite or
countable family of subsets of $\Bbb{R}^n$ with finite metric dimensions bounded by
$\overline{D}_m$, respectively. Let $\{\Delta_m, m\in{\cal M}\}$ be a family of positive
weights satisfying (\ref{Eq-S}). One can build an estimator $\tilde{s}(\bm{X})$ of $s$ such
that, for some universal constant $C$,
\begin{equation}
\Bbb{E}_s\left[\|s-\tilde{s}\|^2\right]\le C\inf_{m\in{\cal M}}
\left\{\sigma^2\max\left\{\overline{D}_m,\Delta_m\right\}+
\inf_{t\in\overline{S}_m}\|s-t\|^2\right\}.
\labe{Eq-pg}
\end{equation}
\end{theorem}
%
%
\begin{theorem}\lab{T-DP}
Let $\bm{X}=(X_1,\ldots,X_n)$ be an i.i.d.\ sample with unknown density $s$ with respect
to some measure  $\nu$ on $(E,{\cal E})$ and $(M,h)$ be the metric space of all such
densities with Hellinger distance. Let $\{\overline{S}_m, m\in{\cal M}\}$ be a finite or
countable family of subsets of $(M,h)$ with finite metric dimensions bounded by
$\overline{D}_m$, respectively. Let $\{\Delta_m, m\in{\cal M}\}$ be a family of positive
weights satisfying (\ref{Eq-S}). One can build an estimator $\tilde{s}(X_1,\ldots,X_n)$ of $s$
such that, for some universal constant $C$,
\begin{equation}
\Bbb{E}_s\left[h^2\left(\tilde{s},s\right)\right]\le C\inf_{m\in{\cal M}}
\left\{n^{-1}\max\left\{\overline{D}_m,\Delta_m\right\}+
\inf_{t\in\overline{S}_m}h^2(s,t)\right\}.
\labe{Eq-pd}
\end{equation}
\end{theorem}
%
\noindent{\it Remark:} The choice of the bound 1 in (\ref{Eq-S}) has nothing canonical and
was simply made for convenience. Any small constant would do since we did not provide
the actual value of  $C$ which depends on the right-hand side of (\ref{Eq-S}).
%
 
\subsubsection{About the complexity of families of models\labs{S2b}}
The only difference between the ideal bound (\ref{Eq-wrb}) and (\ref{Eq-pg}) is the
replacement of $\overline{D}_m$ by $\max\left\{\overline{D}_m,\Delta_m\right\}$ with
weights $\Delta_m$ satisfying (\ref{Eq-S}) and we see, comparing (\ref{Eq-rbd}) and
(\ref{Eq-pd}), that the same difference holds for density estimation. More generally, in a
framework for which an analogue of Proposition~\ref{P-vtest} holds, leading to
(\ref{Eq-gen1}), we proved in Birg\'e (2006) that
\begin{equation}
\Bbb{E}_s\left[d^2(s,\hat{s}_{\hat{m}})\right]\le C
\inf_{m\in{\cal M}}\left\{C_1\inf_{t\in\overline{S}}d^2(s,t)+C_2
\max\left\{\overline{D}_m,\Delta_m\right\}\right\},
\labe{Eq-uopt}
\end{equation}
holds instead of (\ref{Eq-topt}). In all situations, apart from the constant $C$, the loss with
respect to the ideal bound is due to the replacement of $\overline{D}_m$ by
$\max\left\{\overline{D}_m,\Delta_m\right\}$ where the weights $\Delta_m$ satisfy 
(\ref{Eq-S}). If $\Delta_m$ is not much larger than $\overline{D}_m$ for all $m$, we have
almost reached the ideal risk, otherwise not and we can now explain what we mean by the
complexity of a family of models. 

For each positive integer $j$, let us denote by $H(j)$ the cardinality of the set ${\cal M}_j$
of those $m$ such that $j/2\le\overline{D}_m<(j+1)/2$. If $H(j)$ is finite for all $j$, let us 
choose $\Delta_m=(j+1)/2+\log_+(H(j))$ for $m\in{\cal M}_j$ where
$\log_+(x)=\log x$ for $x\ge1$ and $\log_+(0)=0$. Then
\[
\sum_{m\in{\cal M}}\exp[-\Delta_m]=\sum_{j\ge1}\sum_{m\in{\cal M}_j}
\exp[-(j+1)/2-\log_+(H(j))]\le\sum_{i\ge2}\exp[-i/2]<1
\]
and (\ref{Eq-S}) holds. Moreover,
\[
\max\left\{\overline{D}_m,\Delta_m\right\}=\Delta_m\le
2\overline{D}_m[1+j^{-1}\log_+(H(j))]\quad\mbox{for }m\in{\cal M}_j.
\]
If $j^{-1}\log_+[H(j)]$ is uniformly bounded and the bound is not large, then (\ref{Eq-uopt})
and (\ref{Eq-topt}) are comparable and we can consider that the family of models is not
complex. On the other hand, if, for some $j$, $\log[H(j)]$ is substantially larger than
$j$, $\Delta_m$ is substantially larger than $\overline{D}_m$, at least for some $m$,
which may result in a bound (\ref{Eq-uopt}) much larger than (\ref{Eq-topt}). If 
$H(j)=+\infty$ for some $j$, (\ref{Eq-S}) requires that $\Delta_m$ be unbounded for
$m\in{\cal M}_j$, which is even worse. A reasonable measure of the complexity of a
family of models is therefore $\sup_{j\ge1}j^{-1}\log_+[H(j)]$, high complexity of the family 
corresponding to large values of this index.

\subsection{Application 1: variable selection in Gaussian regression\labs{S3}}
Let us now give some concrete illustrations of more or less complex families of models
corresponding to the examples that motivated our investigations about model selection. To
begin with, we consider the situation of Section~\ref{D2a} with a large number $p\le n$ of
potentially influential explanatory variables $Z^j$ and set $\Lambda=\{1;\ldots;p\}$. For
any subset $m$ of $\Lambda$ we define $\overline{S}_m$ as the linear span of the vectors
$Z^j$ for $j\in m$. According to Section~\ref{A3b}, $\overline{S}_m$ has a metric dimension
bounded by
$|m|/2$.

Let us assume that we have ordered the variables according to their supposed relevance,
$Z^1$ being the more relevant. In such a situation it is natural to consider the models spanned
by the $q$ more relevant variables $Z^1,\ldots,Z^q$ for $1\le q\le p$ and therefore to set
${\cal M}={\cal M}_1=\left\{\st\{1;\ldots;q\},1\le q\le p\right\}$. This is not a complex
family of models and  the choice $\Delta_m=|m|$ ensures that (\ref{Eq-S}) holds. It
follows from Theorem~\ref{T-GP} that one can design an estimator $\tilde{s}(\bm{X})$
satisfying
\begin{equation}
\Bbb{E}_s\left[\|s-\tilde{s}\|^2\right]\le C\inf_{m\in{\cal M}_1}
\left\{\sigma^2|m|+\inf_{t\in\overline{S}_m}\|s-t\|^2\right\}.
\labe{Eq-tt5}
\end{equation}
Comparing this with the performance of the m.l.e.\ with respect to each model
$\overline{S}_m$ given by (\ref{Eq-Gauss1}), we see that, apart from the  constant $C$, we
recover the performance of the best model in the family.

This simple approach has, nevertheless, some drawbacks. First, we have to order the
explanatory variables which is often not easy. Then the result is really bad if we make a
serious mistake in this ordering. Imagine, for instance, that $s$ only depends on four
highly influential variables so that if the variables had been ordered correctly, the best model,
i.e.\ the one minimizing $\sigma^2|m|+\inf_{t\in\overline{S}_m}\|s-t\|^2$, would be
$\overline{S}_{\{1;2.3;4\}}$ and the corresponding risk $4\sigma^2$. If one of these four
very influential variables has been neglected and appears in the sequence with a high index
$l$, it may happen that, because of this wrong ordering, the best model becomes
$\overline{S}_{\{1;\ldots;l\}}$ leading to the much higher risk $\sigma^2l$.

In order to avoid the difficulties connected with variables ordering, one may introduce many
more models, defining ${\cal M}={\cal M}_2$ as the set of all nonvoid subsets $m$ of
$\Lambda$.  Since the number of nonvoid subsets of $\Lambda$ with cardinality $q$ is
$\binom{p}{q}\le p^q/q!$, we may choose $\Delta_m=1+|m|\log p$ to get (\ref{Eq-S}) so
that, by Theorem~\ref{T-GP}, one can find an estimator $\tilde{s}(\bm{X})$ satisfying
\begin{equation}
\Bbb{E}_s\left[\|s-\tilde{s}\|^2\right]\le C\inf_{m\in{\cal M}}
\left\{\sigma^2(1+|m|\log p) +\inf_{t\in\overline{S}_m}\|s-t\|^2\right\}.
\labe{Eq-tt6}
\end{equation}
With this method, we avoid the problems connected with variables ordering and may even
introduce more explanatory variables than observations ($p>n$), hoping that with so many
variables at disposal, one can find a small subset $m$ of them that provides an accurate model
for $s$. There is a price to pay for that! We now have a complex family of models when $p$ is
large resulting in values of $\Delta_m$ which are much larger than $|m|$ and we pay the
extra factor $\log p$ in our risk bounds.

One can actually cumulate the advantages of the two approaches by mixing the two families
in the following way. We first order the $p$ variables as we did at the beginning, giving the
smallest indices to the variables we believe are more influential and set again ${\cal M}={\cal
M}_2$. We then fix $\Delta_m=|m|+1/2$ for $m\in{\cal M}_1$ and $\Delta_m=
1+|m|\log p$ for $m\in{\cal M}\setminus{\cal M}_1$ so that (\ref{Eq-S}) still
holds.   Theorem~\ref{T-GP} shows that
\begin{eqnarray*}
\Bbb{E}_s\left[\|s-\tilde{s}\|^2\right]&\le&
C\min\left[\inf_{m\in{\cal M}\setminus{\cal M}_1}\left\{\sigma^2(1+|m|\log p)+
\inf_{t\in\overline{S}_m}\|s-t\|^2\right\}\right.\\&&\mbox{}\qquad\quad\;\;
;\left.\inf_{m\in{\cal M}_1}
\left\{\sigma^2(|m|+1/2)+\inf_{t\in\overline{S}_m}\|s-t\|^2\right\}\right].
\end{eqnarray*}
If our ordering of the variables is right, the best $m$ belongs to ${\cal M}_1$ and we get an
analogue of (\ref{Eq-tt5}). If not, we lose a factor $\log p$ from the risk of the best model as in
(\ref{Eq-tt6}).

\subsection{Application 2: histograms and density estimation\labs{S4}}

\subsubsection{Problems connected with the use of the $\Bbb{L}_2$-distance in
density estimation\labs{S4a}}
Let us now come back to density estimation with histograms. In Section~\ref{H2} we used the
$\Bbb{L}_2$-distance to measure the distortion between $s$ and its estimator. This is certainly
the most popular and more widely studied measure of distortion for density estimation but it
actually has some serious drawbacks as shown by Devroye and Gy\"orfi (1985). For histograms
it results in risk bounds (\ref{Eq-his7}) depending on $\|s\|_\infty$ for irregular partitions,
which are not of the form
\[
R(\hat{s}_{m},s)\le C\left[\|s-s_{m}\|^2+n^{-1}|m|\right],
\]
for some universal constant $C$, independent of $s$, $n$ and the partition $m$. It is
actually impossible to get an analogue of Theorem~\ref{T-Dmodel} where the
$\Bbb{L}_2$-distance would replace the Hellinger distance, as shown by the following
proposition motivated by  Theorem~2.1 of Rigollet and Tsybakov (2005).
Indeed, if such a theorem were true, we could apply it to the model $\overline{S}$ provided
by this proposition and conclude that the corresponding estimator $\hat{s}_{\overline{S}}$
would satisfy the analogue of (\ref{Eq-rbd}) leading to the uniform risk bound
\[
\Bbb{E}_s\left[\|\hat{s}_{\overline{S}}-s\|^2\right]\le CD/(2n),\quad\mbox{for all }s\in
\overline{S}
\]
and some universal constant $C$, therefore independent of $L$. This would clearly contradict
(\ref{Eq-Z2}) below for large enough values of $L$.
%
\begin{proposition}\lab{P-L2}
For each $L>0$ and each integer $D$ with $1\le D\le 3n$, one can find a finite set
$\overline{S}$ of densities with the following properties:

i) it is a subset of some $D$-dimensional affine subspace of $\Bbb{L}_2([0,1],dx)$ with a metric
dimension bounded by $D/2$;

ii) $\sup_{s\in\overline{S}}\|s\|_\infty\le L+1$;

iii) for any estimator $\hat{s}(X_1,\ldots,X_n)$ belonging to $\Bbb{L}_2([0,1],dx)$ and based on
an i.i.d.\ sample with density $s\in\overline{S}$,
\begin{equation}
\sup_{s\in\overline{S}}\Bbb{E}_s\left[\|\hat{s}-s\|^2\right]>0.0139DLn^{-1}.
\labe{Eq-Z2}
\end{equation}
\end{proposition}
%
\noindent{\em Proof:} Let us set  $a=D/(4n)\le3/4$, define $\theta$ by
$(1-\theta)/\theta=4nL/D$ and introduce the functions $f(x)=\1_{[0,1[}(x)$ and
$g(x)=-a\1_{[0,(1-\theta)/D]}+a(1-\theta)\theta^{-1}\1_{](1-\theta)/D],1/D[}$. Then
$\int_0^{1/D}g(x)\,dx=0$, $\sup_xg(x)=L$, $\inf_xg(x)=-a\ge-3/4$ and
\begin{equation}
\|g\|^2=\int_0^{1/D}\!g^2(x)\,dx=a^2\frac{1-\theta}{D}
\left[1+(1-\theta)\theta^{-1}\right]=\frac{a^2(1-\theta)}{\theta D}=\frac{L}{4n}.
\labe{Eq-Z1}
\end{equation}
It follows that $\|f-(f+g)\|^2=L/(4n)$. Moreover 
\begin{eqnarray}
h^2(f,f+g)&=&\frac{1}{2}\int_0^{1/D}\left[1-\sqrt{1+g(x)}\right]^2dx\nonumber\\&=&
\frac{1}{2}\int_0^{1/D}\left[2+g(x)-2\sqrt{1+g(x)}\right]dx\nonumber\;\;=\;\;\frac{1}{D}
-\int_0^{1/D}\sqrt{1+g(x)}\,dx\\&=&\nonumber D^{-1}\left[1-(1-\theta)\sqrt{1-a}-\theta
\sqrt{1+a(1-\theta)\theta^{-1}}\right]\\&\le&D^{-1}\left[1-\sqrt{1-a}\right]\;\;\le\;\;
D^{-1}(2a/3)\;\;=\;\;(6n)^{-1},
\labe{Eq-Z6}
\end{eqnarray}
since $a\le3/4$. Let us now set, for $1\le j\le D$, $g_j(x)=g\left(x-D^{-1}(j-1)\right)$, so that
these $D$ translates of $g$ have disjoint supports and $g_1=g$. Let ${\cal D}=\{0;1\}^D$ with
the distance $\Delta$ given by $\Delta(\delta,\delta')=\sum_{j=1}^D|\delta_j-\delta'_j|$.
For each $\delta\in{\cal D}$ we consider the density $s_\delta(x)=f(x)+\sum_{j=1}^D
\delta_jg_j(x)$ and set $\overline{S}=\{s_\delta,\delta\in{\cal D}\}$. Clearly 
$\|s_\delta\|_\infty\le L+1$ for all $\delta\in{\cal D}$ and it follows from (\ref{Eq-Z1})
that 
\begin{equation}
\|s_\delta-s_{\delta'}\|^2=\sum_{j=1}^D(\delta_j-\delta'_j)^2\int_0^{1/D}\!g_j^2(x)\,dx
=\frac{L}{4n}\sum_{j=1}^D(\delta_j-\delta'_j)^2=\frac{L}{4n}\Delta(\delta,\delta').
\labe{Eq-isom}
\end{equation}
Moreover, since $\overline{S}$ is a subset of some $D$-dimensional affine subspace of
$\Bbb{L}_2([0,1],dx)$, it follows from the arguments used in the proof of
Proposition~\ref{P-mod1} that its metric dimension is bounded by $D/2$.

Defining $P_\delta$ by $dP_\delta/dx=s_\delta$, we derive from (\ref{Eq-Z6}) that
$h^2(P_\delta,P_{\delta'})\le(6n)^{-1}$, hence $\rho(P_\delta,P_{\delta'})\ge\overline{\rho}=
1-(6n)^{-1}$, for  each pair $(\delta,\delta')\in{\cal D}^2$ such that
$\Delta(\delta,\delta')=1$. We may then apply Assouad's Lemma below to conclude from
(\ref{Eq-isom}) that, whatever the estimator $\hat{\delta}$ with values in ${\cal D}$,
\[
\sup_{\delta\in{\cal D}}\Bbb{E}_s\left[\|s_{\hat{\delta}}-s_\delta\|^2\right]=\frac{L}{4n}
\sup_{\delta\in{\cal D}}\Bbb{E}_s\left[\Delta\left(\hat{\delta},\delta\right)\right]
\ge\frac{L}{4n}\frac{D}{2}\left[1-\sqrt{1-\left[1-(6n)^{-1}\right]^{2n}}\right].
\]
Let $\hat{s}$ be any density estimator based on $X_1,\ldots,X_n$ and set
$\hat{\delta}(X_1,\ldots,X_n)$ to satisfy  $\|\hat{s}-s_{\hat{\delta}}\|=
\inf_{\delta\in{\cal D}}\|\hat{s}-s_\delta\|$ so that, whatever $\delta\in{\cal D}$,
$\|s_{\hat{\delta}}-s_\delta\|\le2\|\hat{s}-s_\delta\|$.  We derive from our last bound that
\[
\sup_{\delta\in{\cal D}}\Bbb{E}_s\left[\|\hat{s}-s_\delta\|^2\right]\ge\frac{1}{4}
\sup_{\delta\in{\cal D}}\Bbb{E}_s\left[\|s_{\hat{\delta}}-s_\delta\|^2\right]
\ge\frac{LD}{32n}\left[1-\sqrt{1-\left[1-(6n)^{-1}\right]^{2n}}\right].
\]
We conclude by observing that $\left[1-(6n)^{-1}\right]^{2n}$ is increasing with $n$, hence
$\ge25/36$.\cqfd
%
\begin{lemma} [Assouad, 1983]\lab{L-Assouad}
Let $\left\{P_\delta,\delta\in{\cal D}\right\}$ be a family of distributions indexed by ${\cal
D}=\{0;1\}^D$ and $X_1,\ldots,X_n$ an i.i.d.\ sample from a distribution in the family.
Assume that $\rho(P_\delta,P_{\delta'})\ge\bar{\rho}$ for each pair
$(\delta,\delta')\in{\cal D}^2$ such that $\Delta(\delta,\delta')=1$. Then for any estimator
$\hat{\delta}(X_1,\ldots,X_n)$ with values in ${\cal D}$,
\begin{equation}
\sup_{\delta\in{\cal D}}\Bbb{E}_\delta\left[\Delta\left(\hat{\delta}(X_1,\ldots,X_n),
\delta\right)\right]\ge \frac{D}{2}
\left[1-\sqrt{1-\bar{\rho}^{2n}}\right]\ge\frac{D\bar{\rho}^{2n}}{4},
\labe{Eq-Assou}
\end{equation}
where $\Bbb{E}_\delta$ denotes the expectation when the $X_i$ have the distribution
$P_\delta$.
\end{lemma}
%
\noindent{\em Proof:} 
Let us set $P_\delta^n$ for the joint distribution of the $X_i$ with individual distribution
$P_\delta$ and consider some measure $\mu$ which dominates the probabilities
$P_\delta^n$ for $\delta\in{\cal D}$. First note that the left-hand side of (\ref{Eq-Assou}) is at
least as large as the average risk
\[
R_B=2^{-D}\sum_{\delta\in{\cal D}}\Bbb{E}_\delta\left[\Delta\left(\hat{\delta},
\delta\right)\right]=2^{-D}\sum_{\delta\in{\cal D}}\int\sum_{k=1}^D
\left|\hat{\delta}_k-\delta_k\right|dP_\delta^n.
\]
Then, setting $Q_k^j=2^{-D+1}\sum_{\{\delta\in{\cal D}\,|\,\delta_k=j\}}P_\delta^n$ with
$j=0$ or 1, we get
\begin{eqnarray*}
R_B&=&2^{-D}\sum_{k=1}^D\left(\sum_{\{\delta\in{\cal D}\,|\,\delta_k=0\}}\int
\hat{\delta}_k\,dP_\delta^n+\sum_{\{\delta\in{\cal D}\,|\,\delta_k=1\}}\int\left(1-
\hat{\delta}_k\right)dP_\delta^n\right)\\&=&\frac{1}{2}\sum_{k=1}^D\left(\int
\hat{\delta}_k\frac{dQ_k^0}{d\mu}\,d\mu+\int\left(1-\hat{\delta}_k\right)
\frac{dQ_k^1}{d\mu}\,d\mu\right)\\&\ge&\frac{1}{2}\sum_{k=1}^D
\int\inf\left\{\frac{dQ_k^0}{d\mu};\frac{dQ_k^1}{d\mu}\right\}d\mu.
\end{eqnarray*}
Since $\inf\{x;y\}$ is a concave function of the pair $(x,y)$, it follows that
\[
\inf\left\{\frac{dQ_k^0}{d\mu};\frac{dQ_k^1}{d\mu}\right\}\ge2^{-D+1}
\sum_{(\delta,\delta')\in{\cal D}_k}
\inf\left\{\frac{dP_\delta^n}{d\mu};\frac{dP_{\delta'}^n}{d\mu}\right\},
\]
with ${\cal D}_k=\{(\delta,\delta')\,|\,\delta_k=0,\delta'_k=1,\delta_j=
\delta'_j\mbox{ for }j\ne k\}$, hence
\[
R_B\ge\frac{1}{2}\sum_{k=1}^D2^{-D+1}\sum_{(\delta,\delta')\in{\cal D}_k}\int
\inf\left\{\frac{dP_\delta^n}{d\mu};\frac{dP_{\delta'}^n}{d\mu}\right\}d\mu.
\]
We now use (\ref{Eq-ropi}) to conclude that
\[
R_B\ge\frac{1}{2}\sum_{k=1}^D2^{-D+1}\sum_{(\delta,\delta')\in{\cal D}_k}\left[1-
\sqrt{1-\rho^2\left(P_\delta^n,P_{\delta'}^n\right)}\right].
\]
By assumption, $\rho(P_\delta,P_{\delta'})\ge\bar{\rho}$ for $(\delta,\delta')\in{\cal D}_k$,
hence $\rho^2\left(P_\delta^n,P_{\delta'}^n\right)=\rho^{2n}(P_\delta,P_{\delta'})\ge
\bar{\rho}^{2n}$. The conclusion follows.\cqfd
%
\subsubsection{Partition selection for histograms\labs{S4c}}
If we use the Hellinger distance instead of the $\Bbb{L}_2$-distance to evaluate the risk of
histograms, we can improve (\ref{Eq-his7}), getting a universal bound which does not
involve $\|s\|_\infty$. We recall that $S_m$ is the set of densities which are constant on the
elements of the partition $m$ as defined in (\ref{Eq-Sm}).
%
\begin{theorem}\lab{T-histo}
Let $s$ be some density with respect to the Lebesgue measure on $[0,1]$, $X_1,\ldots,
X_n$ be an $n$-sample from the corresponding distribution and $m=\{I_0,\ldots,I_D\}$ be a
partition of $[0,1]$ into intervals $I_j$ with respective lengths $|I_j|$. Let $\hat{s}_m$ be
the histogram estimator based on this partition and given by
\[
\hat{s}_m(x)=\sum_{j=0}^D\left[\frac{1}{n|I_j|}\sum_{i=1}^n\1_{I_j}(X_i)\right]
\1_{I_j}(x).
\]
The Hellinger risk of $\hat{s}_m$ is bounded by
\begin{equation}
\Bbb{E}_s\left[h^2(s,\hat{s}_m)\right]\le2\inf_{t\in S_m}h^2(s,t)+D/(2n).
\labe{Eq-B1}
\end{equation}
\end{theorem}
%
\noindent{\em Proof:} 
It is shown in Birg\'e and Rozenholc (2006) that
\[ 
\Bbb{E}_s\left[h^2(s,\hat{s}_m)\right]\le h^2(s,s_m)+\frac{D}{2n}\quad\mbox{with }
s_m=\sum_{j=0}^D\left[\frac{1}{|I_j|}\int_{I_j}s(x)\,dx\right]\1_{I_j}. 
\] 
Let $f$ be the $\Bbb{L}_2$-orthogonal projection of $\sqrt{s}$ onto the linear span $V_m$ of
$\1_{I_0},\ldots,\1_{I_D}$. Then 
\[
f=\sum_{j=0}^D\left[\frac{1}{|I_j|}\int_{I_j}\sqrt{s(x)}\,dx\right]\1_{I_j}\quad\mbox{and}\quad
\left\|f-\sqrt{s}\right\|^2\le2h^2(s,t)\quad\mbox{for all }t\in V_m.
\]
Setting $s_m=\sum_{j=0}^Da_j\1_{I_j}$ and $f=\sum_{j=0}^Db_j\1_{I_j}$, we get from
Jensen's Inequality that $b_j\le\sqrt{a_j}$. It follows that
\[
h^2(s,s_m)=1-\sum_{j=0}^D\int_{I_j}\sqrt{a_js(x)}\,dx=1-\sum_{j=0}^D\sqrt{a_j}b_j|I_j|
\le1-\sum_{j=0}^Db_j^2|I_j|,
\]
while
\[
\left\|f-\sqrt{s}\right\|^2=1+\sum_{j=0}^D\int_{I_j}b_j^2\,dx-
2\sum_{j=0}^D\int_{I_j}b_j\sqrt{s(x)}\,dx=1-\sum_{j=0}^Db_j^2|I_j|.
\]
Hence 
\begin{equation}
h^2(s,s_m)\le\left\|f-\sqrt{s}\right\|^2\le2\inf_{t\in S_m}h^2(s,t).\cqfd
\labe{Eq-Hell8}
\end{equation}
%
If, in particular, $\sqrt{s}$ is H\"older continuous and satisfies (\ref{Eq-Hold}), we derive as in
Section~\ref{H3c} that one can find a regular partition $m$, depending on $L$ and $\beta$,
such that,
\begin{equation}
\Bbb{E}_s\left[h^2(s,\hat{s}_m)\right]\le
\max\left\{(5/2)\left(Ln^{-\beta}\right)^{2/(2\beta+1)};n^{-1}\right\}.
\labe{Eq-rhold}
\end{equation}
Then, a useful remark is as follows. If we have at disposal a sample $X_1,\ldots,X_{2n}$ of size
$2n$ and a family ${\cal M}$ of partitions of $[0,1]$, one can use the first half of the sample to
build the corresponding histograms $\hat{s}_m(X_1,\ldots,X_n)$ and use the second half of the
sample to select one estimator in the family. For this, we merely have to apply
Theorem~\ref{T-DP} to the sample $X_{n+1},\ldots,X_{2n}$ conditionally on $X_1,\ldots,X_n$.
Conditionally on $X_1,\ldots,X_n$, each histogram $\hat{s}_m$ is simply a density which can
be considered as a model $\overline{S}_m$ containing only one point, hence with a finite
metric dimension bounded by $1/2$. Let $\{\Delta_m,m\in{\cal M}\}$ be a family of weights 
satisfying 
\begin{equation}
\sum_{m\in{\cal M}}\exp\left[-\Delta_m\right]\le1\qquad\mbox{and}\qquad\Delta_m
\ge1\quad\mbox{for all }m.
\labe{Eq-S'}
\end{equation}
We derive from Theorem~\ref{T-DP} applied to the models $\overline{S}_m=\{\hat{s}_m\}$
that there exists an estimator $\tilde{s}(X_1,\ldots,X_{2n})$ such that
\[
\Bbb{E}_s\left[\left.h^2\left(\tilde{s},s\right)\right|X_1,\ldots,X_n\right]\le
C\inf_{m\in{\cal M}}\left\{n^{-1}\Delta_m+
h^2\left(\st s,\hat{s}_m(X_1,\ldots,X_n)\right)\right\}.
\]
Integrating with respect to $X_1,\ldots,X_n$ and using (\ref{Eq-B1}) finally leads to
\begin{eqnarray}
\Bbb{E}_s\left[h^2\left(\tilde{s},s\right)\right]&\le&\nonumber
C\inf_{m\in{\cal M}}\left\{n^{-1}\Delta_m+2\inf_{t\in S_m}h^2(s,t)
+(|m|-1)/(2n)\right\}\\&\le&C'\inf_{m\in{\cal M}}
\left\{n^{-1}\max\{|m|,\Delta_m\}+\inf_{t\in S_m}h^2(s,t)\right\}.
\labe{Eq-hisf}
\end{eqnarray}

\subsubsection{A straightforward application of partition selection\labs{S4c}}
To give a concrete application of this result, let us introduce some special classes of partitions.
For any finite partition $m=\{I_0,\ldots,I_D\}$ into intervals, we denote by $A_m$ the set
$\{y_0<\ldots<y_{D+1}\}$, $y_0=0, y_{D+1}=1$ of endpoints of the intervals $I_j$.
Introducing, for $k\ge1$, the set ${\cal J}_k$ of dyadic numbers $\{j2^{-k},0\le j\le2^k\}$,
we denote by ${\cal M}_{D,k}$, for $1\le D<2^k$, the set of those partitions $m$ which satisfy
\[
|m|=D+1;\quad A_m\in{\cal J}_k\quad\mbox{and}\quad A_m\not\in{\cal J}_{k-1}.
\]
Denoting by $m_0$ the trivial partition with one element $[0,1]$, we define ${\cal M}$ by
\[
{\cal M}=\{m_0\}\bigcup\left(\bigcup_{k\ge1}\bigcup_{1\le D<2^k}{\cal M}_{D,k}\right).
\]
The partitions in ${\cal M}$ are dense in the set of finite partitions into intervals in the following
sense: given any such partition $m$, an element $t$ in $S_m$, as defined by  (\ref{Eq-Sm}),  and
$\varepsilon>0$, we can find $m'\in{\cal M}$ and $t'\in S_{m'}$ such that
$h(t,t')\le\varepsilon$. This means that the approximation properties of $\bigcup_{m\in{\cal
M}}S_m$ are the same as those of all possible histograms. Since $|{\cal M}_{D,k}|\le
\binom{2^k-1}{D}\le2^{kD}$, if we set $\Delta_m=\Delta^0_m=[(k+1)(D+1)+1]\log2$ for
$m\in{\cal M}_{D,k}$ and $\Delta_{m_0}=1$, we get
\[
\sum_{k\ge1}\sum_{1\le D<2^k}\sum_{{\cal M}_{D,k}}e^{-\Delta^0_m}\le\sum_{k\ge1}
\sum_{1\le
D<2^k}2^{-k-D-2}\le\frac{1}{4}\sum_{k\ge1}2^{-k}\sum_{D\ge1}2^{-D}=\frac{1}{4}.
\]
It follows that (\ref{Eq-S'}) holds so that by (\ref{Eq-hisf}), one can find an estimator 
$\tilde{s}(X_1,\ldots,X_{2n})$ which satisfies
\begin{equation}
\Bbb{E}_s\left[h^2\left(\tilde{s},s\right)\right]\le C\inf_{k\ge1}\inf_{1\le D<2^k}
\inf_{m\in{\cal M}_{D,k}}\left\{\frac{kD}{n}+\inf_{t\in S_m}h^2(s,t)\right\}.
\labe{Eq-ex5}
\end{equation}
If, in the right-hand side of (\ref{Eq-ex5}), we set $m$ to be the regular partition with $2^k$
elements, which belongs to ${\cal M}_{D,2^k-1}$, we get a bound of the form
\[
\Bbb{E}_s\left[h^2\left(\tilde{s},s\right)\right]\le
C\left[k2^kn^{-1}+\inf_{t\in S_m}h^2(s,t)\right].
\]
For densities $s$ with $\sqrt{s}$ satisfying (\ref{Eq-Hold}), we get 
\[
\Bbb{E}_s\left[h^2\left(\tilde{s},s\right)\right]\le C'\inf_{k\ge1}
\left\{k2^kn^{-1}+L^22^{-2k\beta}\right\},
\]
but an optimization with respect to $k$ does not allow to recover the bound (\ref{Eq-rhold})
because of an extra factor $\log\left(nL^2\right)$. This factor is connected with the
complexity  of the families ${\cal M}_{D,k}$ which forces us to fix $\Delta_m$ much larger
than $|m|=D+1$ for most elements of ${\cal M}_{D,k}$ when $k$ is large. Most, but not all! It 
is in particular easy to modify the value of $\Delta_m$ for the regular partitions without
violating  (\ref{Eq-S'}). If $m_k$ denotes the regular partition with $2^k$ elements and ${\cal
M}_R$ the set of such partitions, we may choose $\Delta_{m_k}=|m_k|$ instead of
$\Delta_m^0$ so that 
\[
\sum_{m\in{\cal M}_R}e^{-\Delta_m}=\sum_{k\ge0} e^{-2^k}<0.522
\]
and (\ref{Eq-S'}) still holds.
It is easy to check that, with this new choice of the weights for the regular partitions, we improve
the estimation for those densities such that $\sqrt{s}$ is H\"older continuous. In particular, if
$\sqrt{s}$ satisfies (\ref{Eq-Hold}) for some unknown values of $L$ and $\beta$,
\begin{equation}
\Bbb{E}_s\left[h^2\left(\tilde{s},s\right)\right]\le C
\max\left\{\left(Ln^{-\beta}\right)^{2/(2\beta+1)};n^{-1}\right\},
\labe{Eq-ex6}
\end{equation}
which is comparable to (\ref{Eq-rhold}) although $L$ and $\beta$ are unknown, the only loss
being at the level of the constant $C$.

\subsubsection{Introducing more sophisticated Approximation Theory\labs{S4d}}
The consequences of the previous modification of the weights for partitions in ${\cal M}_R$ is a
simple illustration of the use  of elementary Approximation Theory to improve the estimation of
smooth densities. One can actually do much better with the use of more sophisticated
Approximation Theory. In a milestone paper, Birman and Solomjak (1967) introduced a family
${\cal M}_T$ of partitions of the cube $[0,1]^k$ which are such that piecewise constant (and more
generally piecewise polynomials) based on the partitions in the family have excellent
approximation properties with respect to functions in Sobolev spaces (and functions of bounded
variation when $k=1$). Moreover, Birman and Solomjak provide a control on the number of such
partitions with a given cardinality. For the case $k=1$ which is the one we deal with here, the
number of elements $m$ of ${\cal M}_T$ with $|m|=D$ is bounded $4^D$ which allows us to
set $\Delta_m=2D$ for those partitions.

The algorithm leading to the construction of the partitions in ${\cal M}_T$, which is called an 
``adaptive approximation algorithm", is also described in Section~3.3 of DeVore (1998) and it
works as follows. We choose a positive threshold $\varepsilon$ and some non-negative
functional $J(f,I)$ depending on the function $f$ to be  approximated and the interval $I$.
Roughly speaking, the functional measures the quality of approximation of $f$ by a piecewise
constant (or more generally a piecewise polynomial) function on $I$. At step one, the algorithm
starts with the trivial partition $m^1=m_0$ with one single interval. At step $j$ it provides a
partition $m^j$ into $j$ intervals and it checks whether $\sup_{I\in m^j}J(f,I)\le\varepsilon$ or
not. If this is the case, the algorithm stops, if not we choose one of the intervals $I$ for which the
criterion $J(f,I)\le\varepsilon$ is violated and divide it into two interval of equal length to derive
$m^{j+1}$. Then we iterate the procedure. For the functions $f$ of interest, which satisfy some
smoothness condition related to the functional $J$, the procedure necessarily stops at some
stage, leading to a final partition $m$. Let ${\cal M}_T$ be the set of all the partitions that can be
obtained in this way. Then ${\cal M}_R\subset{\cal M}_T$.  Building a partition $m$ in ${\cal
M}_T$ is actually equivalent to growing a complete binary tree for which the initial interval
$[0,1]$ corresponds to the root of the tree, each node of the tree to an interval and each split of an
interval to adding two sons to a terminal node of the tree, the partition $m$ being in one-to-one
correspondance to the set of terminal nodes of the tree. When viewed as a tree algorithm, this
construction is similar to the CART algorithm of  Breiman, Friedman, Olshen and Stone (1984).
The analysis of CART from the model selection point of view that we explain here has been made
by Gey and N\'ed\'elec (2005).

It follows from the correspondence between the partitions in ${\cal M}_T$ and the complete
binary trees that the number of elements $m$ of ${\cal M}_T$ such that $|m|=j+1$,
$j\in\Bbb{N}$, is equal to the number of complete binary trees with $j+1$ terminal nodes which
is given by the Catalan numbers $(j+1)^{-1}\binom{2j}{j}$. Setting $\Delta^1_m=2|m|$ for
$m\in{\cal M}_T$ and using $\binom{2j}{j}\le4^j$ which follows from Stirling's expansion, we
derive that
\[
\sum_{m\in{\cal M}_T}e^{-\Delta^1_m}\le\sum_{j\ge0}
\frac{e^{-2(j+1)}}{j+1}
\binom{2j}{j}\le\sum_{j\ge0}\frac{4^je^{-2(j+1)}}{j+1}<\frac{1}{4}.
\]
It follows that (\ref{Eq-S'}) holds if we set $\Delta_m=\Delta^1_m$ for $m\in{\cal M_T}$
and $\Delta_m=\Delta^0_m$ for $m\in{\cal M}\setminus{\cal M}_T$ and we then
derive from (\ref{Eq-hisf}) that not only (\ref{Eq-ex5}) still holds but also
\[
\Bbb{E}_s\left[h^2\left(\tilde{s},s\right)\right]\le C\inf_{m\in{\cal M}_T}
\left\{n^{-1}|m|+\inf_{t\in S_m}h^2(s,t)\right\},
\]
which is indeed a substantial improvement over (\ref{Eq-ex5}). In particular, since ${\cal
M}_T$ contains ${\cal M}_R$, (\ref{Eq-ex6}) still holds when $\sqrt{s}$ is H\"olderian, but the
introduction of the much larger class ${\cal M}_T$ leads to a much more powerful result which
follows from the approximation properties of functions in $V_m$ given by (\ref{Eq-Vm}) with
$m\in{\cal M}_T$. We refer the reader to the book by DeVore and Lorentz (1993) for the precise
definitions of Besov spaces and semi-norms and the variation $\Var^*$ in the following
theorem.
%
\begin{theorem}\lab{T-Besappro}
Let ${\cal M}_T$ be the set of partitions $m$ of $[0,1]$ previously defined. For any $p>0$,
$\alpha$ with $1>\alpha>(1/p-1/2)_+$, any positive integer $j$ and any function $t$ 
belonging to the Besov space $B^\alpha_{p,\infty}([0,1])$ with Besov semi-norm
$|t|_{B^\alpha_{p,\infty}}$, one can find some $m\in{\cal M}_T$ with $|m|=j$ and some
$t'\in V_m$ such that
\begin{equation}
\|t-t'\|_2\le C(\alpha,p)|t|_{B^\alpha_{p,\infty}}j^{-\alpha},
\label{Eq-DeVo}
\end{equation}
where $\|\cdot\|_2$ denotes the $\Bbb{L}_2(dx)$-norm on $[0,1]$.

If  $t$ is a function of bounded variation on $[0,1]$, there exists $m\in{\cal M}_T$ with
$|m|=j$  and $t'\in V_m$ such that  $\|t-t'\|_2\le C'\Var^*(t)j^{-1}$.
\end{theorem}
The bound (\ref{Eq-DeVo}) is given in DeVore and Yu (1990). The proof for the bounded
variation case has been kindly communicated to the author by Ron DeVore. 

Applying the previous theorem to $t=\sqrt{s}$, we may always choose for $t'$ the projection of
$\sqrt{s}$ onto $V_m$ and it follows from (\ref{Eq-Hell8}) that the result still holds with
$t'=\sqrt{s_m}$. In particular, if $\sqrt{s}\in B^\alpha_{p,\infty}([0,1])$, then for a
suitable $m$ with $|m|=j$, $h^2(s,s_m)\le C(\alpha,p)|t|^2_{B^\alpha_{p,\infty}}
j^{-2\alpha}$. Putting this into (\ref{Eq-hisf}) with $\Delta_m=2j$ and optimizing with respect
to $j$ shows that
\[
\Bbb{E}_s\left[h^2\left(\tilde{s},s\right)\right]\le C\max\left\{\left(
|t|_{B^\alpha_{p,\infty}}n^{-\alpha}\right)^{2/(2\alpha+1)};n^{-1}\right\}
\quad\mbox{if }\sqrt{s}\in B^\alpha_{p,\infty}([0,1]).
\]
Similarly, we can show that
\[
\Bbb{E}_s\left[h^2\left(\tilde{s},s\right)\right]\le C\max\left\{\left(
\Var^*\left(\sqrt{s}\right)/n\right)^{2/3};n^{-1}\right\}
\quad\mbox{if $\sqrt{s}$ has a bounded variation}.
\]

\subsection{Model choice and Approximation Theory\labs{S5}}
In any statistical framework for which we can prove a risk bound of the form (\ref{Eq-uopt})
provided that (\ref{Eq-S}) holds, the technical problem of model selection can be considered
as being solved but the question of how to choose the family of models to which we shall apply
the procedure remains. There is no general recipe to make such a choice without any ``a priori"
information on $s$. If we have some information about the true $s$ or at least we suspect that it
may have some specific properties, or if we wish that some particular $s$ should be accurately
estimated, we should choose our family of models in such a way that the right-hand side of
(\ref{Eq-uopt}) be as small as possible for the $s$ of interest. Finding models of low dimension
with good approximation properties for some specific functions $s$ is one purpose of
Approximation Theory. One should therefore base our choice of suitable families of models on
Approximation Theory, which accounts for the numerous connections between modern Statistics
and Approximation Theory. 

We may also have the choice between several families of models with different approximation
properties and complexity levels. Typically, the more complex families have better
approximation properties but we have to pay a price for the complexity. A good example is the
alternative regular versus irregular partitions for histograms. As shown in the previous sections,
it is possible to mix families with different approximation and complexity properties by playing
with the weights $\Delta_m$. In particular, it is important that as many models as possible, and
particularly those with good approximation properties with respect to functions of greater
interest, do satisfy $\Delta_m\le c|m|$ for some fixed constant $c$. The introduction of the
family of models $\{S_m, m\in{\cal M}_T\}$ in Section~\ref{S4d} illustrates this fact. These
models, which have especially good approximation properties with respect to a large class of
Besov spaces, form a much richer class than those soleley based on regular partitions.
Nevertheless, the number of such models with dimension $D$ remains bounded by
$\exp[c'D]$, which allows to fix $\Delta_m$ of the order of $D$ for these models. By 
(\ref{Eq-uopt}), this implies that, when we use such a family of models, the performance of the
estimator based on model selection is almost (up to constants) as good as the performance of the
estimator based on the best individual model.

A detailed analysis of the problems of model choice is given in Section~4.1 of Birg\'e and
Massart (2001) which also provides additional information about the relationship between
model selection and Approximation Theory. Further results in this direction are to be found in
Barron, Birg\'e and Massart (1999). It follows from these presentations that all results in
Approximation Theory that describe precisely the approximation properties of some particular
classes of finite dimensional models are of special interest for the statistical applications we
have in mind. Statistics has been using various approximation methods and we would like to
emphasize here two main trends. One is based on approximation of functions by piecewise
polynomials (or similar functions like splines), some major references here being Birman and
Solomjak (1967) and the book by DeVore and Lorentz (1993). The statistical methods
based on this approach to approximation lead to estimators which are generalizations of
histograms, the selection procedure handling the choice of the partition (and also, possibly, the
degree of the polynomials). Another trend is based on the expansion of functions on suitable
bases, formerly the trigonometric basis, more recently bases derived from a multiresolution
analysis (wavelet bases and the like). The related estimators are based on the estimation of the
coefficients in the expansion and the selection chooses the finite set of coefficients to be kept in
the expansion of the final estimator. Statistical procedures based on wavelet thresholding are
of this type. Theorem~\ref{T-Besappro} based on DeVore and Yu (1990) provides a set of
partitions which are relevant for approximation of functions in Besov spaces. A parallel result
by Birg\'e and Massart (2000) applies to the second approach, providing a family of subsets of
coefficients to keep in order to get similar approximation properties.  A  good overview of
nonlinear approximation based on wavelets or piecewise polynomials with many useful
references is to be found in DeVore (1998).

The use of metric entropy or dimensional arguments in Statistics is not new. The first general
results connecting the metric dimension of the parameter set to the performance of estimators
are given by Le Cam (1973 and 1975) and statistical applications of the classical entropy results
by Kolmogorov and Tikhomirov (1961) are developed in Birg\'e (1983). An up to date
presentation with extensions to model selection following ideas by Barron and Cover (1991) is
in Birg\'e (2006). There is also a huge amount of empirical process literature based on entropy
arguments with statistical applications. Many illustrations and references are to be found in
van der Vaart and Wellner (1996), van der Vaart (1998), van de Geer (2000) and Massart (2006).
More generally, connexions between estimation and Approximation Theory, in particular via
wavelet thresholding, have been developed in many papers. Besides the authors' works
already cited, a short selection with further references is as follows: DeVore, Kerkyacharian,
Picard and Temlyakov (2004), Donoho and Johnstone (1994, 1995, 1996 and 1998), Donoho,
Johnstone, Kerkyacharian and Picard (1995, 1996 and 1997), Kerkyacharian and Picard (1992
and 2000) and  Johnstone (1999).\vspace{2mm}\\
\noindent{\bf Acknowledgements} Many thanks to Yannick Baraud, Richard Gill, Piet
Groeneboom, Pascal Massart and Aad van der Vaart for their comments on earlier versions of
this paper.\vspace{5mm}

\noindent{\large{\bf References:}}
{\small

ASSOUAD, P. (1983). Deux remarques sur l'estimation. {\it C. R. Acad. Sc. Paris S\'er. I 
Math.}   {\bf296}, 1021-1024. 

BARRON, A.R., BIRG\'E, L. and MASSART, P. (1999). Risk bounds for model selection via 
penalization.  {\it Probab.\ Theory Relat.\ Fields} {\bf 113}, 301-415.

BARRON, A.R. and COVER, T.M. (1991). Minimum complexity density estimation. {\it 
IEEE Transactions on Information Theory}  {\bf37}, 1034-1054. 

BIRG\'E, L.  (1983). Approximation dans les espaces  m\'etriques   et th\'eorie de 
l'estimation. {\it Z. Wahrscheinlichkeitstheorie Verw. Geb.} {\bf65}, 181-237.  

BIRG\'E, L.  (1984). Stabilit\'e et instabilit\'e du risque minimax pour des variables 
ind\'e\-pen\-dantes \'equidistri\-bu\'ees. {\it Ann.  Inst. H. Poincar\'e Sect. B}  {\bf{20}}, 
201-223.  

BIRG\'E, L.  (1986). On estimating a density using Hellinger distance and some other 
strange facts.  {\it Probab. Theory Relat. Fields} {\bf71}, 271-291.

BIRG\'E, L.  (2004). Model selection for Gaussian regression with random design.
{\it Bernoulli} {\bf 10}, 1039 -1051.

BIRG\'E, L.  (2006). Model selection via testing : an alternative to (penalized) maximum
likelihood estimators. To appear in {\it Ann.\ Inst. Henri Poincar\'e}.

BIRG\'E, L. and MASSART, P. (1993). Rates of convergence for minimum contrast 
estimators. {\it Probab. Th. Rel. Fields}   {\bf 97}, 113-150. 

BIRG\'E, L. and MASSART, P. (1997). From model selection to adaptive estimation. In {\it 
Festschrift for Lucien Le Cam: Research Papers in Probability and Statistics} (D. Pollard, E.
Torgersen and G. Yang, eds.), 55-87. Springer-Verlag, New York. 

BIRG\'E, L. and MASSART, P. (1998). Minimum contrast estimators on sieves: exponential 
bounds and rates of convergence. {\it Bernoulli} {\bf 4}, 329-375.

BIRG\'E, L. and MASSART, P.  (2000). An adaptive compression algorithm in Besov
spaces.  {\it Constructive Approximation} {\bf 16} 1-36.

BIRG\'E, L. and MASSART, P. (2001). Gaussian model selection. {\it J. Eur. Math. Soc.}
{\bf 3}, 203-268.

BIRG\'E, L. and ROZENHOLC, Y. (2006). How many bins should be put in a regular
histogram. {\it ESAIM-PS} {\bf 10}, 24-45.

BIRMAN, M.S. and SOLOMJAK, M.Z. (1967). Piecewise-polynomial approximation of
functions of the classes $W_{p}$.   {\it Mat.  Sbornik} {\bf 73}, 295-317.  

BREIMAN, L., FRIEDMAN, J.H., OLSHEN, R.A.\ and STONE, C.J. (1984). {\it Classification
and Regression Trees}. Wadsworth, Belmont.

CENCOV, N.N.  (1962). Evaluation of an unknown distribution density from observations. 
{\it Soviet Math.}  {\bf3}, 1559-1562.

CRAM\'ER, H. (1946). {\it Mathematical Methods of Statistics}. Princeton University
Press, Princeton.

DeVORE, R.A. (1998). Nonlinear Approximation. {\it Acta Numerica} {\bf 7}, 51-150.

DeVORE, R.A., KERKYACHARIAN, G., PICARD, D.\ and TEMLYAKOV, V. (2004).
Mathematical methods for supervised learning. Tech.\ report 0422, IMI, University of South
Carolina, Columbia.

DeVORE, R.A.\ and LORENTZ, G.G. (1993). {\it Constructive Approximation}. 
Springer-Verlag, Berlin. 

DeVORE, R.A.\ and YU,Ê X.M. (1990). Degree of adaptive approximation. {\it Math. Comp.}
{\bf  55}, 625-635.

DEVROYE, L.  and GY\"ORFI, L. (1985). {\it Nonparametric Density Estimation: The $L_1$
View}. John Wiley, New York.

DONOHO, D.L. and JOHNSTONE, I.M. (1994). Ideal spatial adaptation by wavelet 
shrinkage.  {\it Biometrika} {\bf 81}, 425-455.

DONOHO, D.L. and JOHNSTONE, I.M. (1995). Adapting to unknown smoothness via 
wavelet shrinkage. {\it JASA} {\bf 90}, 1200-1224.

DONOHO, D.L. and JOHNSTONE, I.M. (1996). Neo-classical minimax problems,
thresholding and adaptive function estimation. {\it Bernoulli} {\bf 2}, 39-62.

DONOHO, D.L. and JOHNSTONE, I.M. (1998). Minimax estimation via wavelet shrinkage. 
 {\it Ann. Statist.} {\bf 26}, 879-921.

DONOHO, D.L., JOHNSTONE, I.M., KERKYACHARIAN, G. and PICARD, D. (1995). 
Wavelet shrinkage: Asymptopia? {\it J. R. Statist. Soc.} B  {\bf57}, 301-369.  

DONOHO, D.L., JOHNSTONE, I.M., KERKYACHARIAN, G. and PICARD, D. (1996). 
Density estimation by wavelet thresholding {\it Ann. Statist.} {\bf 24}, 508-539.

DONOHO, D.L., JOHNSTONE, I.M., KERKYACHARIAN, G. and PICARD, D. (1997).
Universal near minimaxity of wavelet shrinkage. In {\it  Festschrift for Lucien Le Cam:
Research Papers in Probability and Statistics} (D. Pollard, E. Torgersen and G. Yang, eds.),
183-218. Springer-Verlag, New York. 

FISCHER, R.A. (1922). On the mathematical foundations of theoretical statistics. {\it Philos.
Trans. Royal Soc. London Ser. A} {\bf 222}, 309-368.

FISCHER, R.A. (1925). Theory of statistical estimation. {\it Proc. Cambridge Philos. Soc.}
{\bf 22}, 700-725.

GEY, S. and N\'ED\'ELEC, E. (2005). Model selection for CART regression trees. {\it IEEE
Transactions on Information Theory} {\bf  51}, 658-670.

GRENANDER, U. (1981). {\it Abstract inference}. John Wiley, New York. 

GROENEBOOM, P. and WELLNER, J.A. (1992). {\it Information Bounds and
Nonparametric Maximum Likelihood Estimation}. Birkh\"auser, Basel.

GY\"ORFI, L., KOHLER, M., KRY\.ZAK, A. and WALK, H. (2002). {\it A
Distribution-Free Theory of Nonparametric Regression}. Springer, New York.

H\'AJEK, J. (1970). A characterization of limiting distributions of regular estimates. {\it Z.
Wahrsch. Verw. Gebiete} {\bf 14}, 323-330.

H\'AJEK, J. (1972). Local asymptotic minimax and admissibility in estimation. {\it Proc. Sixth
Berkeley Symp. Math. Statist. Probab.} {\bf 1}, 174-194. Univ. California Press, Berkeley.

IBRAGIMOV, I.A. and HAS'MINSKII, R.Z.  (1981). {\it Statistical Estimation: Asymptotic 
Theory}. Springer-Verlag, New York.

JOHNSTONE, I. (1999). {\it Function Estimation and Gaussian Sequence Models}. Book in
preparation.\\
{\sf http://www-stat.stanford.edu/people/faculty/johnstone/baseb.pdf}

KERKYACHARIAN, G. and PICARD, D. (1992). Density estimation in Besov spaces. {\it 
Statist. and  Probab. Lett.} {\bf 13}, 15-24.

KERKYACHARIAN, G. and PICARD, D. (2000). Thresholding algorithms, maxisets and
well-concentrated bases. {\it Test} {\bf 9}, 283-344. 

KOLMOGOROV, A.N. and TIKHOMIROV, V.M. (1961). $\varepsilon $-entropy and   
$\varepsilon$-capacity of sets in function spaces. {\it  Amer. Math. Soc. Transl. (2)} {\bf17},
277-364.  

Le CAM, L.M. (1953). On some asymptotic properties of maximum likelihood estimates and
related Bayes' estimates. {\it Univ. California Publ. Statist.} {\bf 1}, 277-329. 

Le CAM, L.M. (1970). On the assumptions used to prove asymptotic normality of maximum
likelihood estimates. {\it Ann. Math. Stat.} {\bf  41}, 802-828. 

Le CAM, L.M. (1973). Convergence of estimates under dimensionality restrictions. {\it
Ann. Statist.}  {\bf1} , 38-53. 

Le CAM, L.M. (1975). On local and global properties in the theory of asymptotic normality of
experiments. {\it Stochastic Processes and Related Topics, Vol.\ 1} (M. Puri, ed.), 13-54. 
Academic Press,  New York. 

Le CAM, L.M. and YANG, G.L. (2000). {\it Asymptotics in Statistics: Some Basic Concepts.
Second Edition}. Springer-Verlag, New York.

MASSART, P.  (2006). Conentration Inequalities and Model Selection. In {\it  Lecture on
Probability Theory and Statistics, Ecole d'Et\'e de  Probabilit\'es de Saint-Flour XXXIII -
2003} (J.~Picard, ed.). Lecture Note in Mathematics, Springer-Verlag,  Berlin. 

RIGOLLET, T. and TSYBAKOV, A.B. (2005). Linear and convex aggregation of density
estimators. Technical report, University Paris VI.

SHEN, X. and WONG, W.H. (1994). Convergence rates of sieve estimates. {\it Ann. Statist.} 
{\bf22}, 580-615.  

SILVERMAN, B.W.  (1982). On the estimation of a probability density function by the 
maximum penalized likelihood method. {\it Ann. Statist.} {\bf10}, 795-810.  

van de GEER, S. (1990). Estimating a regression function. {\it Ann. Statist.} {\bf18},  
907-924. 

van de GEER, S. (1993). Hellinger-consistency of certain nonparametric maximum 
likelihood estimates. {\it Ann. Statist.} {\bf21}, 14-44.

van de GEER, S. (1995). The method of sieves and minimum contrast estimators. {\it Math. 
Methods Statist.\ }{\bf 4}, 20-38.

van de GEER, S. (2000). {\it Empirical Processes in $M$-Estimation}. Cambridge
University Press,  Cambridge. 

van der VAART, A.W. (1998). {\it Asymptotic Statistics}. Cambridge University Press, 
Cambridge. 

van der VAART, A.W. (2002). The statistical work of Lucien Le Cam. {\it Ann. Statist.} {\bf30},
631-682.

van der VAART, A.W. and WELLNER, J.A. (1996). {\it Weak Convergence and Empirical
Processes, With Applications to Statistics}. Springer-Verlag,  New York.

WAHBA, G. (1990).  {\it Spline Models for Observational Data}. S.I.A.M., Philadelphia.

WALD, A. (1949). Note on the consistency of the maximum likelihood estimate. {\it Ann. 
Math. Statist.}  {\bf20} , 595-601. 

WONG, W.H. and SHEN, X. (1995). Probability inequalities for likelihood ratios and 
convergence rates of sieve MLEs.  {\it Ann. Statist.}  {\bf23}, 339-362. 
}\vspace{4mm}\\
Lucien BIRG\'E\\ UMR 7599 ``Probabilit\'es et mod\`eles al\'eatoires"\\
Laboratoire de Probabilit\'es, bo\^{\i}te 188\\
Universit\'e Paris VI, 4 Place Jussieu\\
F-75252 Paris Cedex 05\\
France\vspace{2mm}\\
e-mail: LB@CCR.JUSSIEU.FR

\end{document}